\documentclass[a4paper,12pt]{article}
\usepackage{epsfig}
\usepackage{amsmath}
\usepackage{amsfonts}
\usepackage{amssymb}
\usepackage{mathrsfs}

\usepackage{color}

\usepackage{amsfonts}
\usepackage{graphicx}
\usepackage{textcomp}
\usepackage{euscript}
\usepackage{makeidx}

\baselineskip
\usepackage{amssymb}
\usepackage{epsfig}
\usepackage{graphicx}
\newtheorem{Lemma}{Lemma}[section]
\newtheorem{Theorem}[Lemma]{Theorem}

\begin{document}

	\title{\Large On the formulation of Adams-Bashforth scheme with Atangana-Baleanu-Caputo fractional derivative to model chaotic problems\footnotemark}

\author{\small Kolade M. Owolabi\footnotemark \;\;and Abdon Atangana\footnotemark \\
	\small \footnotemark[2]\;\footnotemark[3] Institute for Groundwater Studies, Faculty of Natural and Agricultural Sciences \\
	\small University of the Free State, Bloemfontein 9300, South
	Africa\\
	\small \footnotemark[2]  Department of Mathematical Sciences,  Federal University of Technology, PMB 704, Akure, \\
	\small Ondo State, Nigeria\\}
\date{}
\maketitle
\def\thefootnote{\fnsymbol{footnote}}
\setcounter{footnote}{0} \footnotetext[2]{E-mail addresses:
	kmowolabi@futa.edu.ng (K.M. Owolabi); abdonatangana@yahoo.fr (A. Atangana)}
\noindent

	\begin{abstract}
		\noindent
    Mathematical analysis with numerical application of the newly formulated fractional version of the Adams-Bashforth method using the Atangana-Baleanu derivative which has nonlocal and nonsingular properties is considered in this paper. We adopt the fixed point theory and approximation method to prove the existence and uniqueness of the solution via a general two-component time fractional differential equations. The method is tested with three nonlinear chaotic dynamical system in which the integer-order derivative is modelled with the proposed fractional-order case. Simulation result for different $\alpha$ values in $(0,1]$ is presented.
	\end{abstract}

\section*{Leading paragraph}
A range of chaotic and hyperchaotic processes were modeled with the Atangana-Baleanu fractional derivative which has both nonlocal and nonsingular properties in the sense of Caputo. The standard or local derivative in such systems are replaced in time with the fractional derivative versions. A two-step fractional Adams-Bashforth scheme is proposed to numerically approximate the Atangana-Baleanu operator. We apply the fixed point theory and approximation method to show the existence and uniqueness of solutions when applied to a two-component dynamical system.  Numerical results are given for different instances of fractional power $\alpha$.
	
\section{Introduction}
Fractional calculus is regarded as applicable mathematics. The theories and properties of  these fractional operators have generated a lot of interest, and become an active subjects of study in the last few decades. Recently, scientists, engineers and applied  mathematicians have found the fractional calculus concept useful in various fields: quantitative biology, rheology, electro-chemistry, diffusion process, scattering theory, transport theory, elasticity, probability and potential theory \cite{Aza16, Cap16, Pod99, Spr10}. In-fact, it has been regarded as the future of mathematical modelling of real-life phenomena in groundwater studies, geo-hydrology and fractals \cite{Ata16f, Ata17,Ata17a,Ata17b,Ata18,Ata18a}.

However, many engineers, technologists, mathematicians and other scientists are still not familiar with this area of research possibly because they have not been used to its applications. Thus, while the theory of  fractional calculus has been developed, its  use  is still poorly understood or lagged behind. So, another achievement of this paper  is to encourage researchers to find and apply the new numerical scheme \cite{Ata-Owo,Owo18,Owo18a,Owo18b,Owo18c,Owo18d,Owo18e,Owo18f,Owo18g,Owo18h} with the Atangana-Baleanu (AB) derivative in the Caputo sense, which represents some chaotic phenomena with mathematical equations that can  be treated with the robustness of fractional calculus.

Chaotic processes are best described as nonlinear dynamical problems which are more sensitive to initial data, and dense periodic orbits. In recent years, many chaotic systems  such as Zhou-system, Zhu-system, Li-system, Wei-Yang system,  Sundarapandian-system, Vaidyanathan-system and Pham-system to mention a few. A list of old and new chaotic systems are is well classified in \cite{Aza16, Pet11, Spr10}. Chaos and control dynamics have a lot of application in almost every areas of applied science and engineering \cite{Ata16f, Ata17b, Aza15}. Mostly known applications are biology, chemical reactions, ecology, groundwater treatment, lasers, oscillators and finance \cite{Ata16f, Ata17, Ata16, Gom16, Gom16a, Gom17, Owo16a, Owo17, Owo17a, Owo17}, and references therein.

Specifically speaking, we consider a time-fractional $\alpha-$th order nonlinear system
\begin{equation}\label{gen-system}
\begin{split}
\mathcal{D}_{0,t}^\alpha u_{i}(t)&=f_i(u_1(t),u_2(t),\ldots,u_{n}(t),t)\\
u_{i}(0) &=\Phi_i,\;\;\;i=1,2,\ldots,n,
\end{split}
\end{equation}
where $\alpha\in(0,1]$ determines the order of the differential system, $\Phi_i$ is the initial functions, and nonlinear (coupled) functions  $f_i(u_1,u_2,\ldots,u_{n},t)$ for $i=1,2,\ldots,n$ which account for the local kinetics of the system. The fractional derivative $\mathcal{D}_{0,t}^\alpha u_{i}(t)$ is defined by the Atangana-Baleanu derivative
\begin{equation}
^{ABC}_0\mathcal{D}_{t}^\alpha u_{i}(t)=\frac{1-\alpha}{AB(\alpha)}f_i(u_1,u_2,\ldots,u_{n},t)
+\frac{\alpha}{AB(\alpha)\Gamma(\alpha)}\int_{0}^{t}f_i(u_1,u_2,\ldots,u_{n},\tau)(t-\tau)^{\alpha-1}d\tau.
\end{equation}

The rest of this paper is broken into different sections. A quick tour of some basic properties of fractional differentiation and integration is reviewed in Section 2. By following the work reported in \cite{Ata-Owo}, we present the fractional version of the proposed method in the sense of Caputo in Section 3, and show that the non-integer derivative is well-defined and Lipschitz continuous. Numerical experiment via three chaotic fractional differential equations is considered in Section 4. Finally, we conclude with the last section.

\section{Basic properties}
Here, we quickly have a review of some necessary and useful properties of fractional calculus, based on fractional integration and differentiation for a continuous function $u\in C^n[a,b]$. Let $u:[a,b]\rightarrow \mathbb{R}$ be a function, $\alpha>0$ a non-integer value and $n\in\mathbb{N}$ be such that $\alpha\in(n-1, n)$. We assumed that all the necessary conditions for the fractional operator to be well defined is satisfied by $u$.

The left-right Riemann-Liouville (RL) integrals of fractional-order $\alpha\in(0,1]$ is respectively defined as
$$^{RL}\mathcal{I}_{a,t}^\alpha u(t)=\frac{1}{\Gamma(\alpha)}\int_{a}^{t}(t-\xi)^{\alpha-1}u(\xi)d\xi$$
and
$$^{RL}\mathcal{I}_{t,b}^\alpha u(t)=\frac{1}{\Gamma(\alpha)}\int_{t}^{b}(\xi-t)^{\alpha-1}u(\xi)d\xi.$$

The left-right RL derivatives of fractional-order $\alpha\in(0,]$ is respectively defined as
$$^{RL}\mathcal{D}_{a,t}^\alpha u(t)=\frac{1}{\Gamma(n-\alpha)}\frac{d^n}{dt^n}\int_{a}^{t}(t-\xi)^{n-\alpha-1}u(\xi)d\xi$$
and
$$^{RL}\mathcal{D}_{t,b}^\alpha u(t)=\frac{1}{\Gamma(n-\alpha)}\frac{d^n}{dt^n}\int_{t}^{b}(\xi-t)^{n-\alpha-1}u(\xi)d\xi.$$

Likewise, the left-right Caputo fractional derivatives are defined by
$$^C\mathcal{D}_{a,t}u(t)=\frac{1}{\Gamma(n-\alpha)}\int_{a}^{t}(t-\xi)^{n-\alpha-1}u^{(n)}(\xi)d\xi$$
and
$$^C\mathcal{D}_{t,b}u(t)=\frac{1}{\Gamma(n-\alpha)}\int_{a}^{t}(\xi-t)^{n-\alpha-1}u^{(n)}(\xi)d\xi$$
respectively, where $n>0$ is an integer not less than $\alpha$.
The connections between the left-right RL and the Caputo operators of fractional order $0<\alpha\le 1$ is respectively given by
$$^C\mathcal{D}_{a,t}u(t)={^{RL}\mathcal{D}_{a,t}u(t)}-\sum_{s=0}^{n-1}\frac{u^{(s)(a)}}{\Gamma(s-\alpha+1)}(t-a)^{s-\alpha},$$
and
$$^C\mathcal{D}_{t,b}u(t)={^{RL}\mathcal{D}_{t,b}u(t)}-\sum_{s=0}^{n-1}\frac{u^{(s)(b)}}{\Gamma(s-\alpha+1)}(b-t)^{s-\alpha}.$$
Hence, if
$u(a)=u'(a)=\cdots =u^{(n-1)}(a)=0$, then $^C\mathcal{D}_{a,t}u(t)={^{RL}\mathcal{D}_{a,t}u(t)}$. Similarly, if $u(b)=u'(b)=\cdots =u^{(n-1)}(b)=0$, then $^C\mathcal{D}_{t,b}u(t)={^{RL}\mathcal{D}_{t,b}u(t)}$, and $u^{(s)}$ integrable on both $[a,t]$ and $[t,b]$, see \cite{Kil06, Pod99} for details.

In-line with the Caputo, a derivative with nonsingular and nonlocal properties which was proposed by Atangana and Baleanu in 2016 \cite{Ata16a} as
		\begin{equation}\label{ABC}
		^{ABC}\mathcal{D}_{a,t}^\alpha[u(t)]=\frac{M(\alpha)}{1-\alpha}\int_{a}^{t}u'(\xi)E_\alpha\left[-\alpha\frac{(t-\xi)^\alpha}{1-\alpha}\right]d\xi
		\end{equation}
where $M(\alpha)$ is the usual Caputo-Fabrizio normalized function \cite{Cap15, Cap16}, and $E_\alpha$ is a one-parameter Mittag-Leffler function expressed in terms of power series
\begin{equation}\label{Mittag}
u(z)=E_\alpha(z)=\sum_{s=0}^{\infty}\frac{z^s}{\Gamma(\alpha s+1)}, \;\;\;\alpha>0,\;\;\alpha\in\mathbb{R},\;\;z\in\mathbb{C}.
\end{equation}

\section{Approximation of the Atangana-Baleanu fractional derivative}
We examine the following differential equation
\begin{equation}\label{ab1}
_0^{ABC}\mathsf{D}_t^\alpha u(t)=f(t,u(t)).
\end{equation}
By applying the fundamental theorem, we have
\begin{equation}\label{ab2}
u(t)-u(0)=\frac{1-\alpha}{ABC(\alpha)}f(t,u(t))+\frac{\alpha}{ABC(\alpha)\Gamma(\alpha)}\int_{0}^{t}(t-\tau)^{\alpha-1}f(\tau,u(\tau))d\tau.
\end{equation}
At $t_{n+1}$, we have
$$u(t_{n+1})-u(0)=\frac{1-\alpha}{ABC(\alpha)}f(t_n,u_n)+\frac{\alpha}{ABC(\alpha)\Gamma(\alpha)}\int_{0}^{t_{n+1}}(t_{n+1}-\tau)^{\alpha-1}f(t,u(t))dt $$
and at $t_n$ we have
$$u(t_{n})-u(0)=\frac{1-\alpha}{ABC(\alpha)}f(t_{n-1},u_{n-1})+\frac{\alpha}{ABC(\alpha)\Gamma(\alpha)}\int_{0}^{t_{n}}(t_{n}-\tau)^{\alpha-1}f(t,u(t))dt $$
which on subtraction yields
\begin{eqnarray}
u(t_{n+1})-u(t_n)&=&\frac{1-\alpha}{ABC(\alpha)}\left\{f(t_n,u_n)-f(t_{n-1},u_{n-1})\right\}+\frac{\alpha}{ABC(\alpha)\Gamma(\alpha)}\int_{0}^{t_{n+1}}(t_{n+1}-t)^{\alpha-1} \nonumber\\
&&\times f(t,u(t))dt -  \frac{\alpha}{ABC(\alpha)\Gamma(\alpha)}\int_{0}^{t_{n}}(t_{n}-t)^{\alpha-1}f(t,u(t))dt.
\end{eqnarray}
Therefore,
$$u(t_{n+1})-u(t_n)=\frac{1-\alpha}{ABC(\alpha)}\left\{f(t_n,u_n)-f(t_{n-1},u_{n-1})\right\}+A_{\alpha,1}-A_{\alpha,2}.$$

Next, we consider
$$A_{\alpha,1}=\frac{\alpha}{ABC(\alpha)\Gamma(\alpha)}\int_{0}^{t_{n+1}}(t_{n+1}-t)^{\alpha-1}f(t,u(t))dt$$
Again we use the approximation
\begin{equation}
p(t)=\frac{t-t_{n-1}}{t_n-t_{n-1}}f(t_n,u_n)+ \frac{t-t_{n-1}}{t_{n-1}-t_{n}}f(t_{n-1},u_{n-1})
\end{equation}
thus
\begin{eqnarray}
A_{\alpha,1}&=&\frac{\alpha}{ABC(\alpha)\Gamma(\alpha)}\int_{0}^{t_{n+1}}(t_{n+1}-t)^{\alpha-1}\left\{\frac{t-t_{n-1}}{h}f(t_n,u_n)-\frac{t-t_n}{h}f(t_n,u_n)  \right\}\nonumber\\
&=&\frac{\alpha f(t_n,u_n)}{ABC(\alpha)\Gamma(\alpha)h}\left\{\int_{0}^{t_{n+1}}(t_{n+1}-t)^{\alpha-1}f(t-t_{n-1}) \right\}dt\nonumber\\
&& - \frac{\alpha f(t_{n-1},u_{n-1})}{ABC(\alpha)\Gamma(\alpha)h}\left\{\int_{0}^{t_{n+1}}(t_{n+1}-t)^{\alpha-1}f(t-t_{n-1}) \right\}dt\nonumber\\
&=&\frac{\alpha f(t_n,u_n)}{ABC(\alpha)\Gamma(\alpha)h}\left\{\frac{2h t_{n+1}^\alpha}{\alpha}-\frac{ t^{\alpha+1}_{n+1}}{\alpha+1} \right\} - \frac{\alpha f(t_{n-1},u_{n-1})}{ABC(\alpha)\Gamma(\alpha)h}\left\{\frac{h t_{n+1}^\alpha}{\alpha}-\frac{ t^{\alpha+1}_{n+1}}{\alpha+1} \right\}.
\end{eqnarray}
Similarly, we obtain
\begin{equation}
A_{\alpha,2}=\frac{\alpha f(t_n,u_n)}{ABC(\alpha)\Gamma(\alpha)h} \left\{\frac{h t_{n}^\alpha}{\alpha}-\frac{ t^{\alpha+1}_{n}}{\alpha+1} \right\}-\frac{ f(t_{n-1},u_{n-1})}{ABC(\alpha)\Gamma(\alpha)h}
\end{equation}
thus
\begin{eqnarray}
u(t_{n+1})-u(t_n)&=&\frac{1-\alpha}{ABC(\alpha)} \left\{f(t_n,u_n)-f(t_{n-1},u_{n-1})\right\} +\frac{\alpha f(t_n,u_n)}{ABC(\alpha)\Gamma(\alpha)h}\left\{\frac{2h t_{n+1}^\alpha}{\alpha}-\frac{t^{\alpha+1}_{n+1}}{\alpha+1}\right\}\nonumber\\
&&-\frac{\alpha f(t_{n-1},u_{n-1})}{ABC(\alpha)\Gamma(\alpha)h}\left\{\frac{h t_{n+1}^\alpha}{\alpha}-\frac{t^{\alpha+1}_{n+1}}{\alpha+1}\right\}-
\frac{\alpha f(t_n,u_n)}{ABC(\alpha)\Gamma(\alpha)h}\left\{\frac{h t_{n}^\alpha}{\alpha}-\frac{t^{\alpha+1}_{n}}{\alpha+1}\right\}\nonumber\\
&&+\frac{ f(t_{n-1},u_{n-1})}{ABC(\alpha)\Gamma(\alpha)}t_n^{\alpha+1}.
\end{eqnarray}

\begin{eqnarray}
u(_{n+1})&=&u(t_n) + f(t_n,u_n)\left\{\frac{1-\alpha}{AB(\alpha)}+\frac{\alpha}{ABC(\alpha)h}\left[\frac{2h t_{n+1}^\alpha}{\alpha}-\frac{t^{\alpha+1}_{n+1}}{\alpha+1}  \right]\right.\nonumber\\
&&\left.-\frac{\alpha}{AB(\alpha)\Gamma(\alpha)h}\left[\frac{ht^\alpha_n}{\alpha}-\frac{t^{\alpha+1}_n}{\alpha+1}\right]\right\} +f(t_{n-1},u_{n-1}) \\
&&\times\left\{\frac{\alpha-1}{AB(\alpha)}-\frac{\alpha}{h\Gamma(\alpha)AB(\alpha)}\left[\frac{ht_{n+1}^\alpha}{\alpha}-\frac{t_{n+1}^{\alpha+1}}{\alpha+1}+ \frac{t^{\alpha+1}}{h\Gamma(\alpha)AB(\alpha)} \right]  \right\}.\nonumber
\end{eqnarray}
This equation is known as the two-step fractional Adams-Bashforth scheme for Atangana-Baleanu  derivative in  Caputo sense.

Further, we simplify the last equation by putting $t_n=nh$ and $t_{n+1}=(n+1)h$, and collect the resulting equation in powers of $h$ to yield
\begin{eqnarray}\label{eqn-ABC}
u_{n+1}&=&u_n + f(t_n,u_n)\left\{\frac{1-\alpha}{AB(\alpha)}-\frac{\alpha}{AB(\alpha)\Gamma(\alpha)}h^{\alpha}\left[\frac{2(n+1)^\alpha }{\alpha}-\frac{(n+1)^{\alpha+1}}{\alpha+1}\right]\right.\nonumber\\
&&\left.-\frac{\alpha}{AB(\alpha)\Gamma(\alpha)} h^{\alpha}\left[\frac{n^\alpha}{\alpha}-\frac{n^{\alpha+1}}{\alpha+1}\right]   \right\}+ f(t_{n-1},u_{n-1})\\
&&\times\left\{\frac{\alpha-1}{AB(\alpha)}-\frac{\alpha}{AB(\alpha)\Gamma(\alpha)}h^{\alpha}\left[\frac{(n+1)^\alpha}{\alpha}-\frac{(n+1)^{\alpha+1}}{\alpha+1}+\frac{n^{\alpha+1}}{AB(\alpha)\Gamma(\alpha)h}\right]  \right\}.\nonumber
\end{eqnarray}\\
In the followings we briefly outline the convergence and stability results previously obtained in \cite{Ata-Owo}.

\begin{Theorem}{\bf(Convergence result \cite{Ata-Owo})}
	Let $u(t)$ be a solution of
	$$^{ABC}_0\mathcal{D}_t^\alpha u(t)=f(t,u(t))$$
	with function $f$ continuous and bounded, the  solution of u(t) is defined as
\begin{eqnarray*}
u_{n+1}&=&u_n + f(t_n,u_n)\left\{\frac{1-\alpha}{AB(\alpha)}+\frac{\alpha}{AB(\alpha)h}\left[\frac{2h t_{n+1}^\alpha}{\alpha}-\frac{t^{\alpha+1}_{n+1}}{\alpha+1}  \right]\right.\nonumber\\
&&\left.-\frac{\alpha}{AB(\alpha)\Gamma(\alpha)h}\left[\frac{ht^\alpha_n}{\alpha}-\frac{t^{\alpha+1}_n}{\alpha+1}\right]\right\} +f(t_{n-1},u_{n-1}) \\
&&\times\left\{\frac{\alpha-1}{AB(\alpha)}-\frac{\alpha}{h\Gamma(\alpha)AB(\alpha)}\left[\frac{ht_{n+1}^\alpha}{\alpha}-\frac{t_{n+1}^{\alpha+1}}{\alpha+1}+ \frac{t^{\alpha+1}}{h\Gamma(\alpha)AB(\alpha)} \right]  \right\}+R_\alpha
\end{eqnarray*}
where $\|R_\alpha\|_\infty<M$
\end{Theorem}

\begin{Theorem}{\bf (Condition for stability \cite{Ata-Owo})}
	If function $f$ is Lipschitz continuous, then the stability condition for fractional version of the Adams-Bashforth scheme when applied to approximate the ABC derivative  is obtained if
	$$\left\|f(t_n,u_n)-f(t_{n-1},u_{n-1})\right\|_\infty\rightarrow 0$$ as $n\rightarrow \infty$.
\end{Theorem}

\section{Uniqueness and existence of solution via chaotic process}
In this section, we choose a chaotic system to examine the uniqueness and existence of the new method when applied to approximate the ABC fractional derivative. To start with, we consider the general two component system
\begin{equation}
\begin{split}
_0^{ABC}\mathcal{D}_t^\alpha x(t)&=f(x,y,t),\\
_0^{ABC}\mathcal{D}_t^\alpha y(t)&=g(x,y,t),
\end{split}
\end{equation}
By adopting the fundamental calculus theorem along $x$and $y$ components, we have
\begin{equation}
\begin{split}
x(t)-x(0)&=\frac{1-\alpha}{AB(\alpha)}f(x,y,t)+\frac{\alpha}{AB(\alpha)\Gamma(\alpha)}\int_{0}^{t}(t-\tau)^{\alpha-1}f(x,y,\tau)d\tau\\
y(t)-y(0)&=\frac{1-\alpha}{AB(\alpha)}g(x,y,t)+\frac{\alpha}{AB(\alpha)\Gamma(\alpha)}\int_{0}^{t}(t-\tau)^{\alpha-1}g(x,y,\tau)d\tau
\end{split}
\end{equation}
	
Next, we require to create a compact $G_{a,b}$, which means
\begin{equation}
G_{a,b}=I_a(t_0)\times \mathcal{B}_b(\xi)
\end{equation}
where
$$\xi=\min\{x_0,x_0\}$$
and
$$I_a(t_0)=[t_0-a, t_0+a],\;\;\;\mathcal{B}_0(\xi)=[\xi-b,\xi+b].$$
Let
$$M=\max_{G_{a,b}}\left\{\sup_{G_{a,b}}\|f\|,\;\sup_{G_{a,b}}\|g\|\right\}.$$
By adopting the infinite norm, we get
$$\|\varPhi\|_\infty=\sup_{t\in I_a}||\varPhi(t)|.$$
Also, we create a function, say
$$\Gamma : G_{a,b}\rightarrow G_{a,b}$$
so that
\begin{equation}
\begin{split}
\Gamma x(t)&=x_0+\frac{1-\alpha}{AB(\alpha)}f(x,y,t)+\frac{\alpha}{AB(\alpha)\Gamma(\alpha)}\int_{0}^{t} f(x,y,t)(t-\tau)^{\alpha-1}d\tau\\
\Gamma y(t)&=y_0+\frac{1-\alpha}{AB(\alpha)}g(x,y,t)+\frac{\alpha}{AB(\alpha)\Gamma(\alpha)}\int_{0}^{t} g(x,y,t)(t-\tau)^{\alpha-1}d\tau
\end{split}
\end{equation}

Next to prove that the fractional operator is well-defined, we evaluate the condition for which
\begin{equation}
\begin{split}
	\left\|\Gamma_1 x(t)-x_0 \right\|_\infty&<b,\\
	\left\|\Gamma_2 y(t)-y_0 \right\|_\infty&<b.
\end{split}
\end{equation}

So, beginning with the $x$ component, we get
\begin{equation}
\begin{split}
\left\|\Gamma_1 x(t)-x_0 \right\|_\infty=&\left\|\frac{1-\alpha}{AB(\alpha)}f(x,y,t)+\frac{\alpha}{AB(\alpha)\Gamma(\alpha)}\int_{0}^{t} f(x,y,t)(t-\tau)^{\alpha-1}d\tau \right\|_\infty\\
\le&\frac{1-\alpha}{AB(\alpha)}\|f(x,y,t)\|_\infty+\frac{\alpha}{AB(\alpha)\Gamma(\alpha)}\|f(x,y,t)\int_{0}^{t}(t-\tau)d\tau \\
\le &\frac{(1-\alpha)M}{AB(\alpha)}+\frac{\alpha M}{AB(\alpha)\Gamma(\alpha)}\cdot {a^\alpha}<b
\end{split}
\end{equation}
This implies that
$$a=\left(\frac{b-\frac{(1-\alpha)M}{AB(\alpha)}}{\frac{\alpha M}{AB(\alpha)\Gamma(\alpha)}}\right)^{\frac{1}{\alpha}}.$$



Next, we show that the functions $x(t)$ and $y(t)$ satisfy a Lipschitz condition. Which implies
\begin{equation}
\|\Gamma x_1-\Gamma x_2\|_{\infty}\le K \|x_1-x_2\|_{\infty}
\end{equation}
which implies,
\begin{equation}
\begin{split}
\Gamma(x_1)=&\frac{1-\alpha}{AB(\alpha)}f(x_1,y,t)+\frac{\alpha}{AB(\alpha)\Gamma(\alpha)}\int_{0}^{t}f(x_1,y,\tau)(t-\tau)^{\alpha-1}d\tau,\\
\Gamma(x_2)=&\frac{1-\alpha}{AB(\alpha)}f(x_2,y,t)+\frac{\alpha}{AB(\alpha)\Gamma(\alpha)}\int_{0}^{t}f(x_2,y,\tau)(t-\tau)^{\alpha-1}d\tau,
\end{split}
\end{equation}
so that
\begin{eqnarray}
\begin{split}
\|\Gamma x_1-\Gamma x_2\|_{\infty}=&\frac{1-\alpha}{AB(\alpha)}\|f(x_1,y,t)-f(x_2,y,t)\|_{\infty}+\frac{\alpha}{AB(\alpha)\Gamma(\alpha)}\|f(x_1,y,t)\\
&-f(x_2,y,t)\|_{\infty}\int_{0}^{t}(t-\tau)d\tau\\
\le&\|f(x_1,y,t)-f(x_2,y,t)\|_{\infty}\left(\frac{1-\alpha}{AB(\alpha)}+\frac{\alpha}{AB(\alpha)\Gamma(\alpha)}\cdot\frac{a^\alpha}{\alpha}\right)\\
\le&\|f(x_1,y,t)-f(x_2,y,t)\|_{\infty}\left(\frac{1-\alpha}{AB(\alpha)}+\frac{a^\alpha}{AB(\alpha)\Gamma(\alpha)}\right)
\end{split}
\end{eqnarray}
In other words, if $f$ is Lipschitz continuous with respect to $x$, then
\begin{equation}
\begin{split}
\|\Gamma x_1-\Gamma x_2\|_{\infty}\le& K\|x_1-x_2\|_{\infty}\left\{\frac{1-\alpha}{AB(\alpha)}+\frac{a^\alpha}{AB(\alpha)\Gamma(\alpha)}\right\}\\
\le& L \|x_1-x_2\|_{\infty}.
\end{split}
\end{equation}
Similarly, $g$ is Lipschitz with respect to $y$ if
\begin{equation}
\begin{split}
\|\Gamma y_1-\Gamma y_2\|_{\infty}\le&  L \|y_1-y_2\|_{\infty}.
\end{split}
\end{equation}
The procedure described above can be generalized to multicomponent case of fractional differential equations.

\section{Numerical experiments}
Under this segment, we apply the numerical scheme illustrated in the above section to some chaotic systems. The classical time-derivative is replaced with the Atangana-Baleanu fractional derivative of order $\alpha\in(0,1]$ in the sense of Caputo.

From (\ref{eqn-ABC}), we have
\begin{equation}\label{comp}
u_{n+1}=u_n +\omega_1(n,\alpha,h)f(t_n,u_n)+\omega_2(n,\alpha,h)f(t_{n-1},u_{n-1})
\end{equation}
where
\begin{eqnarray}
\omega_1(n,\alpha,h)&=&\left\{\frac{1-\alpha}{AB(\alpha)}-\frac{\alpha}{AB(\alpha)\Gamma(\alpha)}h^{\alpha}\left[\frac{2(n+1)^\alpha }{\alpha}-\frac{(n+1)^{\alpha+1}}{\alpha+1}\right]\right.\nonumber\\
&&\left.-\frac{\alpha}{AB(\alpha)\Gamma(\alpha)} h^{\alpha}\left[\frac{n^\alpha}{\alpha}-\frac{n^{\alpha+1}}{\alpha+1}\right]   \right\}\nonumber
\end{eqnarray}
and
\begin{eqnarray}
\omega_2(n,\alpha,h)&=&f(t_{n-1},u_{n-1})\left\{\frac{\alpha-1}{AB(\alpha)}\right.\nonumber\\
&&\left.-\frac{\alpha}{AB(\alpha)\Gamma(\alpha)}h^{\alpha}\left[\frac{(n+1)^\alpha}{\alpha}-\frac{(n+1)^{\alpha+1}}{\alpha+1}+\frac{n^{\alpha+1}}{AB(\alpha)\Gamma(\alpha)h}\right]  \right\}.\nonumber
\end{eqnarray}

\subsection{Example 1} The following chaotic system \cite{Aza16} is described by the AB fractional derivative in the sense of Caputo as
\begin{equation}\label{exp1}
\begin{split}
^{ABC}_0\mathcal{D}_t^\alpha x_1(t)=&f_1(x_1,x_2,x_3)=\phi(x_2(t)-x_1(t))+\sigma x_2(t) x_3(t),\\
^{ABC}_0\mathcal{D}_t^\alpha x_2(t)=&f_2(x_1,x_2,x_3)=\varphi x_1(t)-x_1(t) x_3(t),\\
^{ABC}_0\mathcal{D}_t^\alpha x_3(t)=&f_3(x_1,x_2,x_3)=x_1(t)x_2(t)-\psi x_3(t)+\delta x_2^2(t),
\end{split}
\end{equation}
where $x_1(t), x_2(t), x_3(t)$ are the densities or states, and $\phi,\varphi,\psi,\sigma,\delta$ are positive parameters.
From (\ref{comp}), we have
\begin{eqnarray}
x_{1,n+1}&=&x_{1,n}+\omega_1(n,\alpha,h)f_1(t_n,x_{1,n},x_{2,n},x_{3,n}) +\omega_2(n,\alpha,h)f_1(t_{n-1},x_{1,n-1},x_{2,n-1},x_{3,n-1})\nonumber\\
x_{2,n+1}&=&x_{2,n}+\omega_1(n,\alpha,h)f_2(t_n,x_{1,n},x_{2,n},x_{3,n}) +\omega_2(n,\alpha,h)f_2(t_{n-1},x_{1,n-1},x_{2,n-1},x_{3,n-1})\nonumber\\
x_{3,n+1}&=&x_{3,n}+\omega_1(n,\alpha,h)f_3(t_n,x_{1,n},x_{2,n},x_{3,n}) +\omega_2(n,\alpha,h)f_3(t_{n-1},x_{1,n-1},x_{2,n-1},x_{3,n-1})\nonumber
\end{eqnarray}

\begin{figure}[!h]
	\begin{center}
	\begin{tabular}{cc}
		\begin{minipage}{200pt}
		\includegraphics[width=200pt]{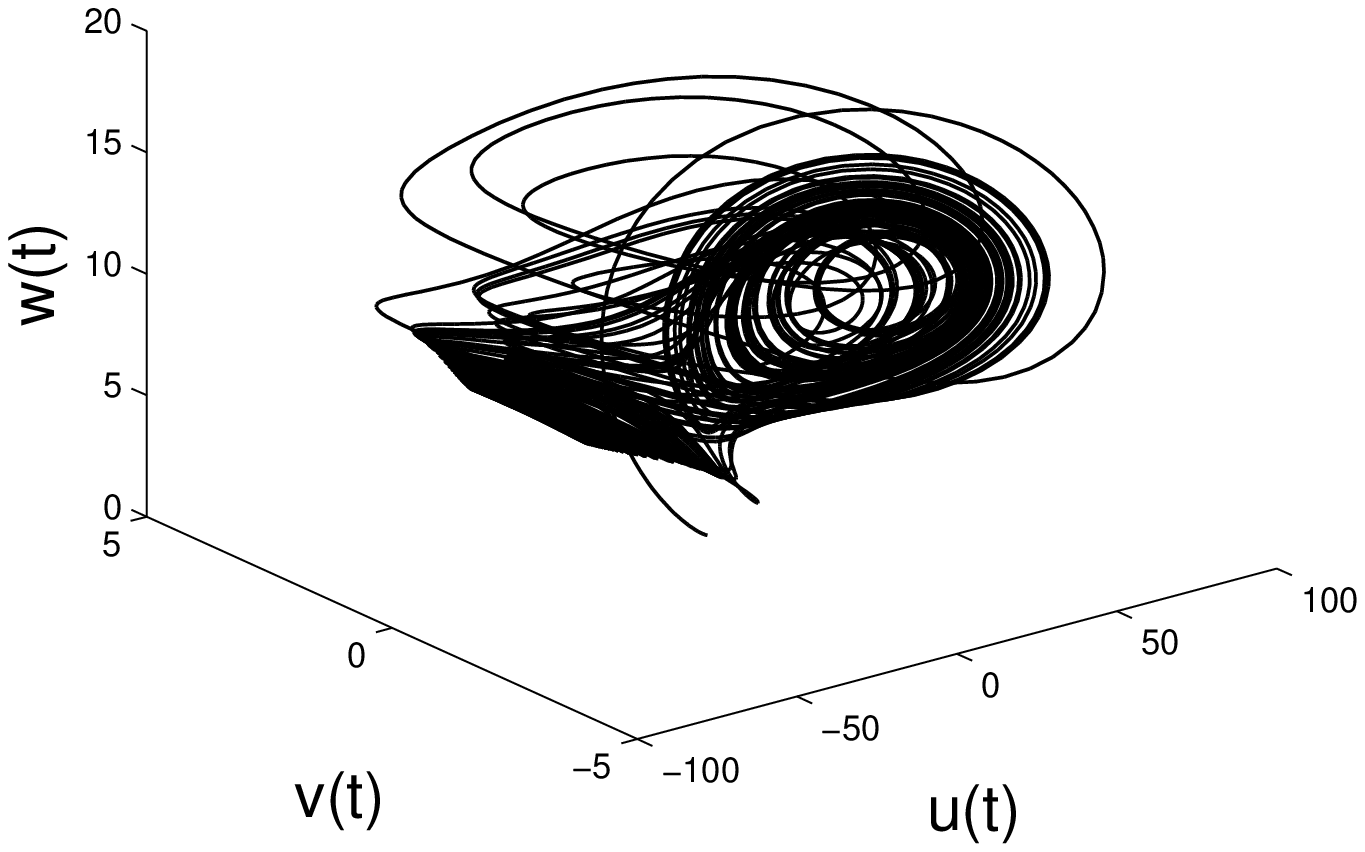}
		\end{minipage}
		\begin{minipage}{200pt}
		\includegraphics[width=200pt]{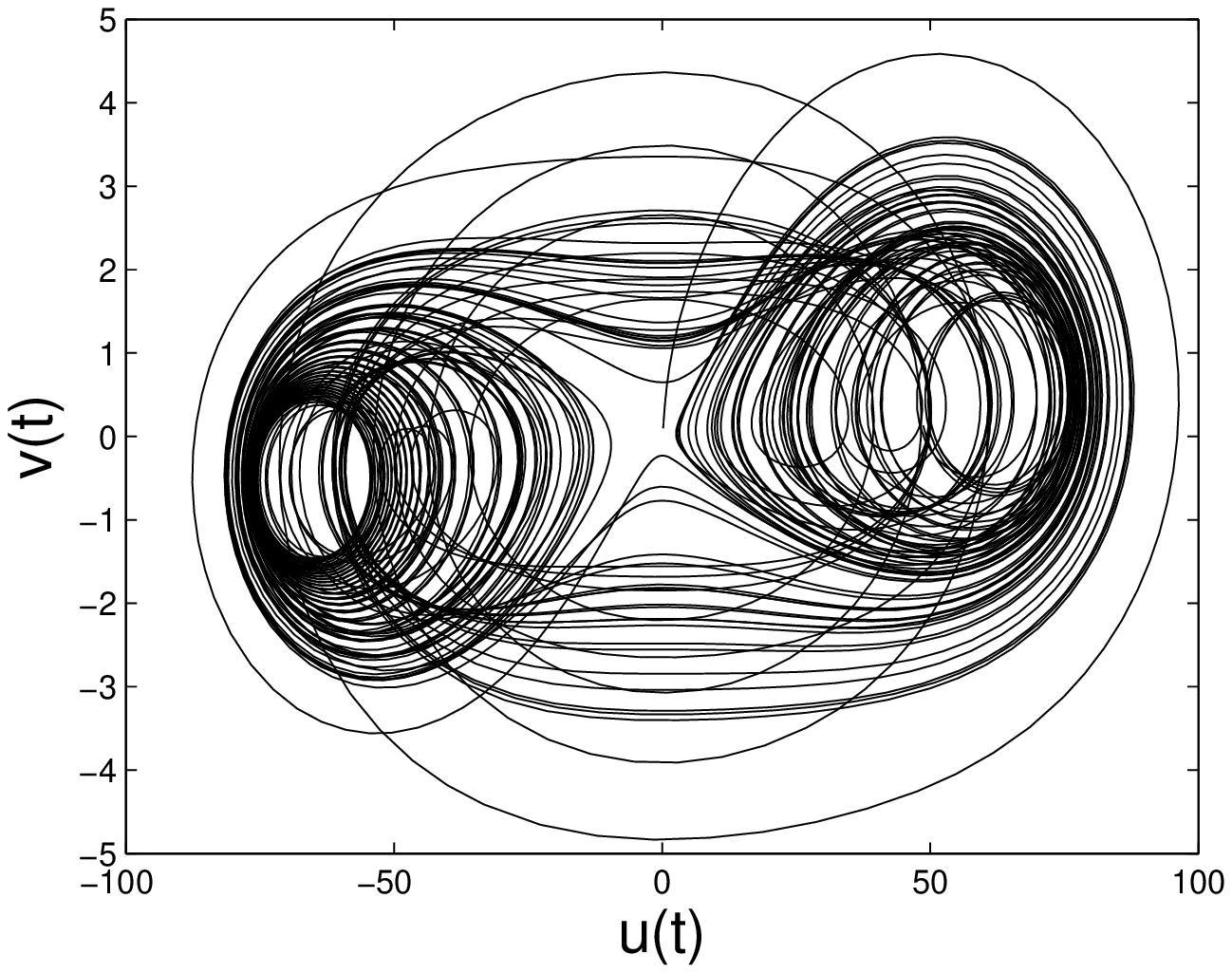}
		\end{minipage}\\
		\begin{minipage}{200pt}
		\includegraphics[width=200pt]{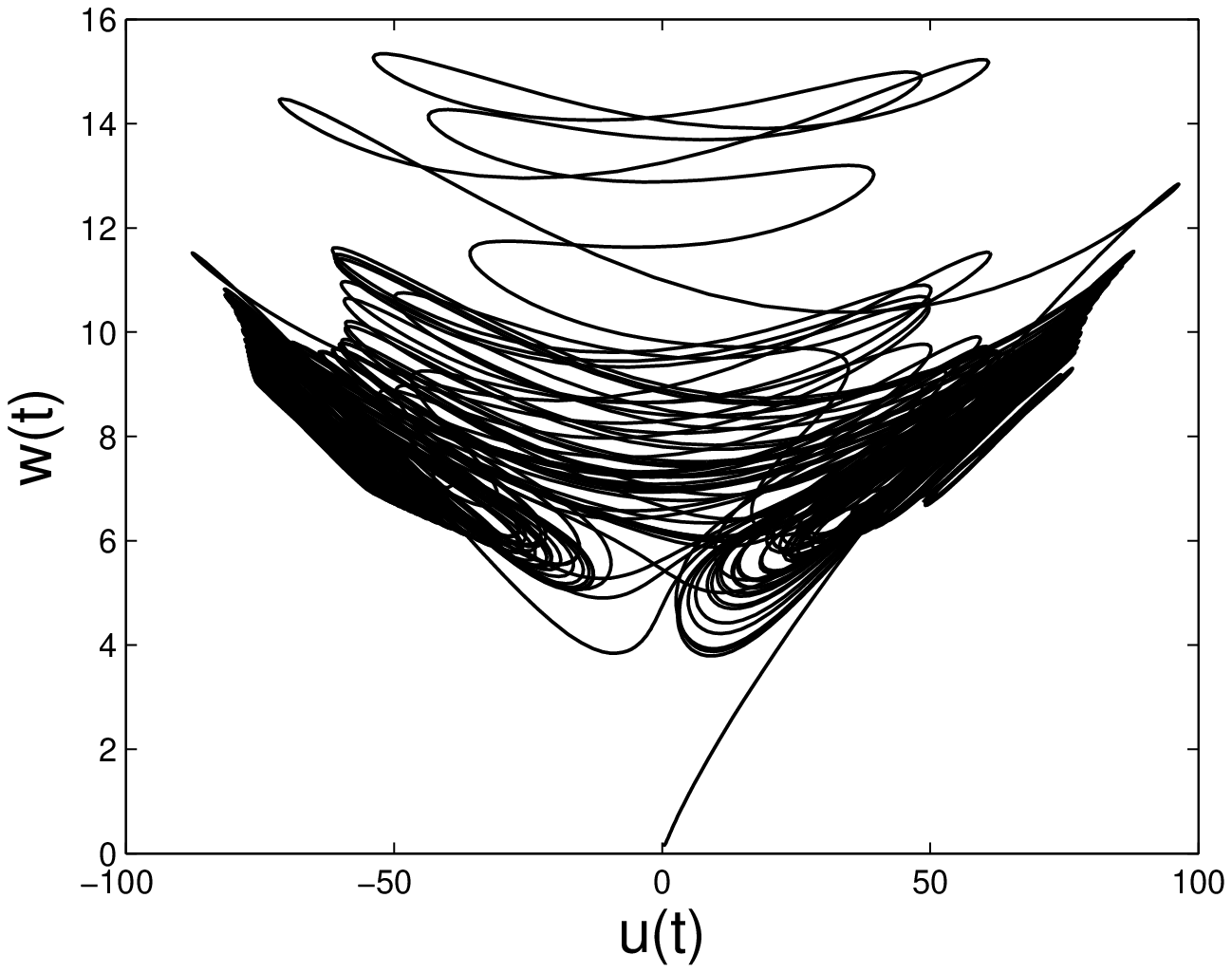}
		\end{minipage}
		\begin{minipage}{200pt}
		\includegraphics[width=200pt]{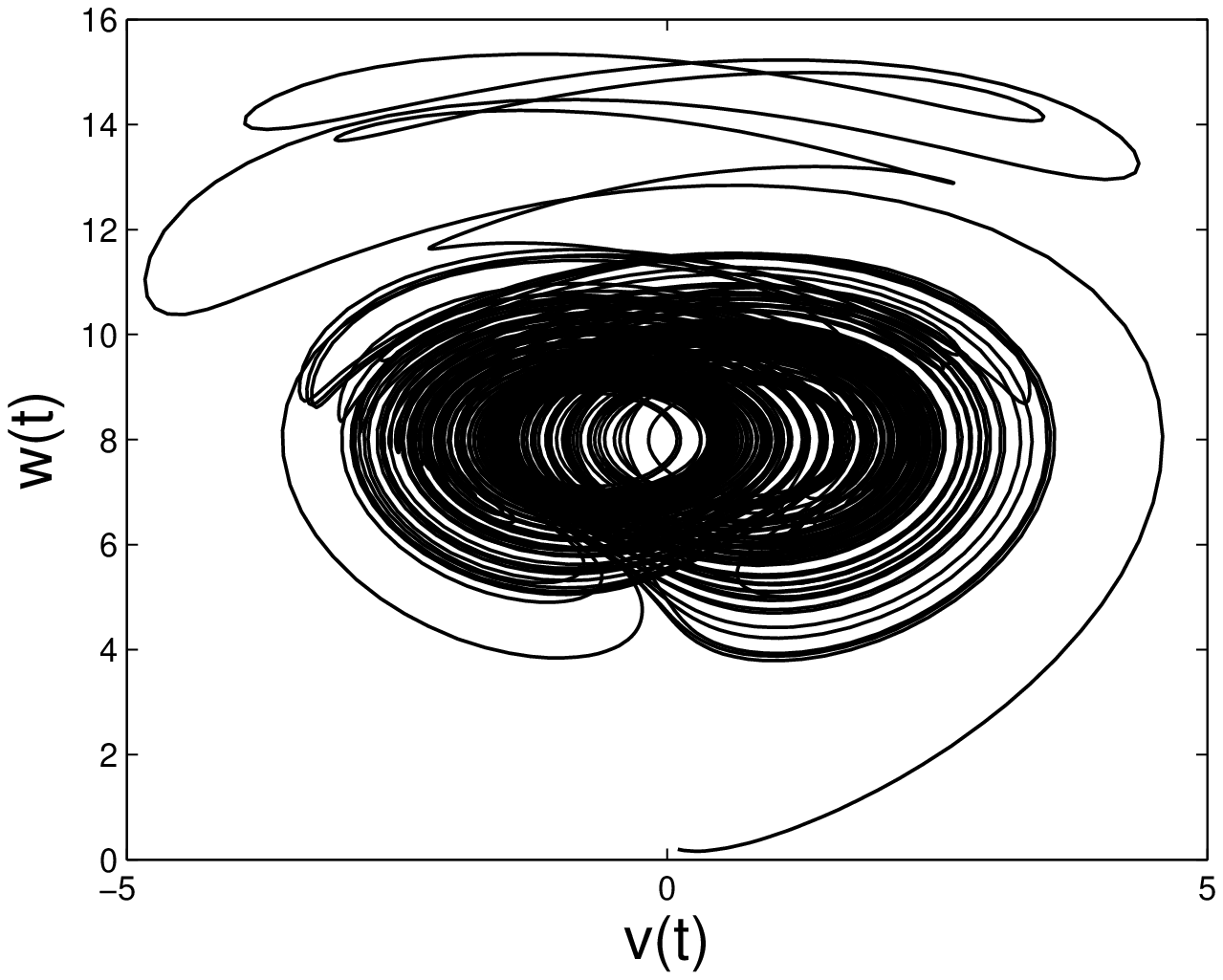}
		\end{minipage}	
	\end{tabular}
    \end{center}
	\caption{Numerical experiment results for chaotic system (\ref{exp1}) obtained at $\alpha=0.75$.}\label{Fig1a}
\end{figure}

\begin{figure}[!h]
	\begin{center}
	\begin{tabular}{cc}
		\begin{minipage}{200pt}
		\includegraphics[width=200pt]{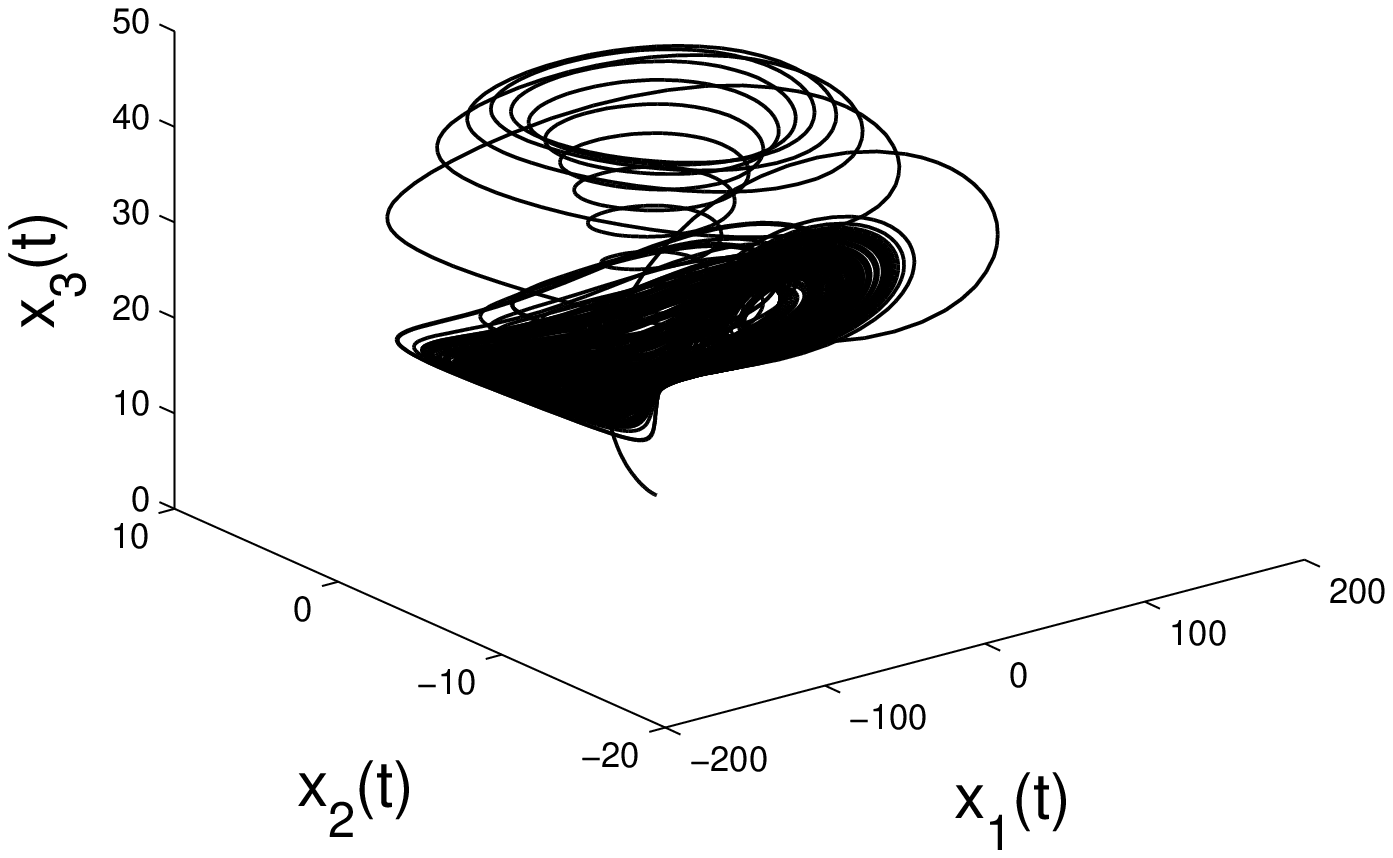}
		\end{minipage}
		\begin{minipage}{200pt}
		\includegraphics[width=200pt]{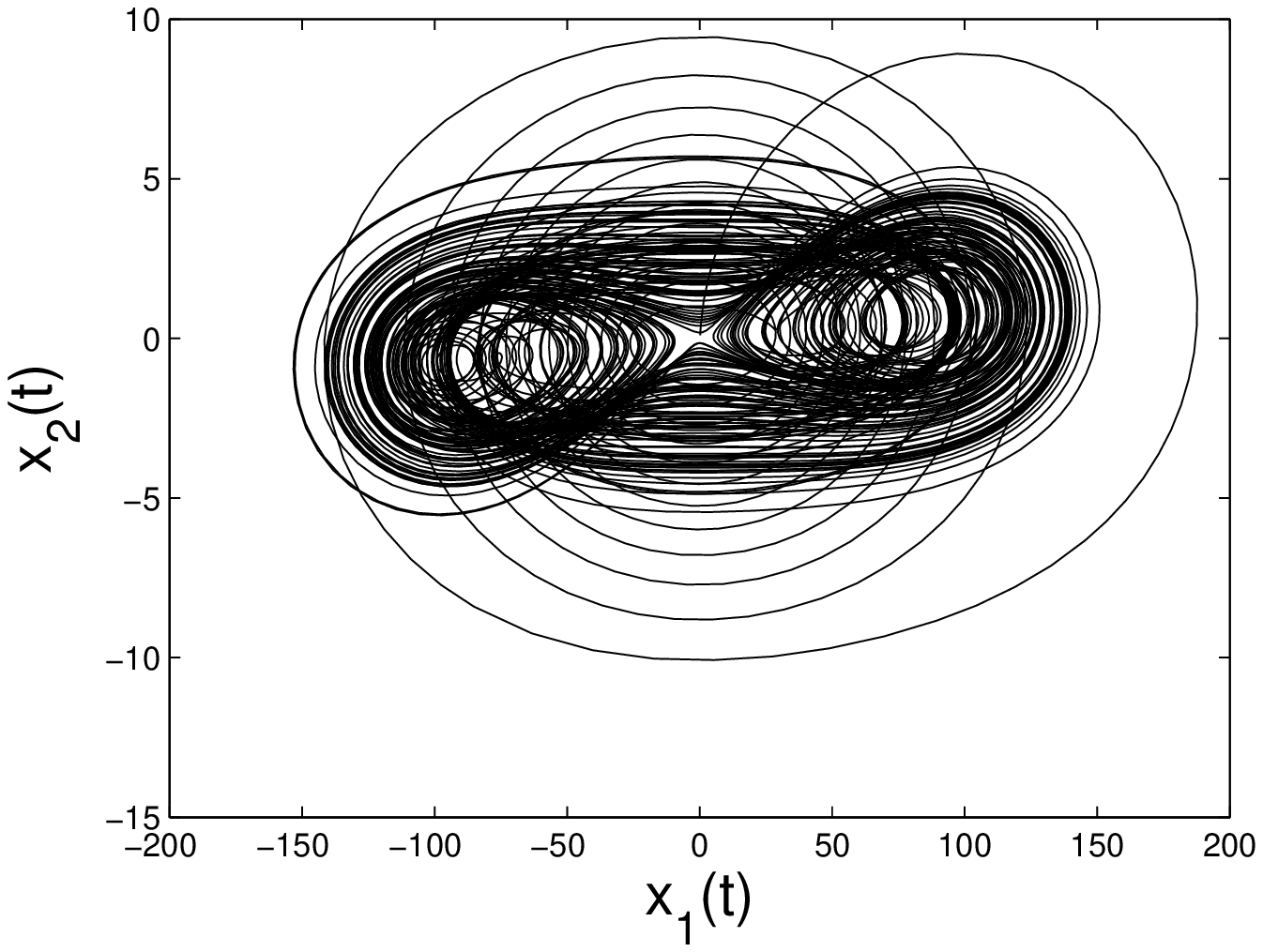}
		\end{minipage}\\
		\begin{minipage}{200pt}
		\includegraphics[width=200pt]{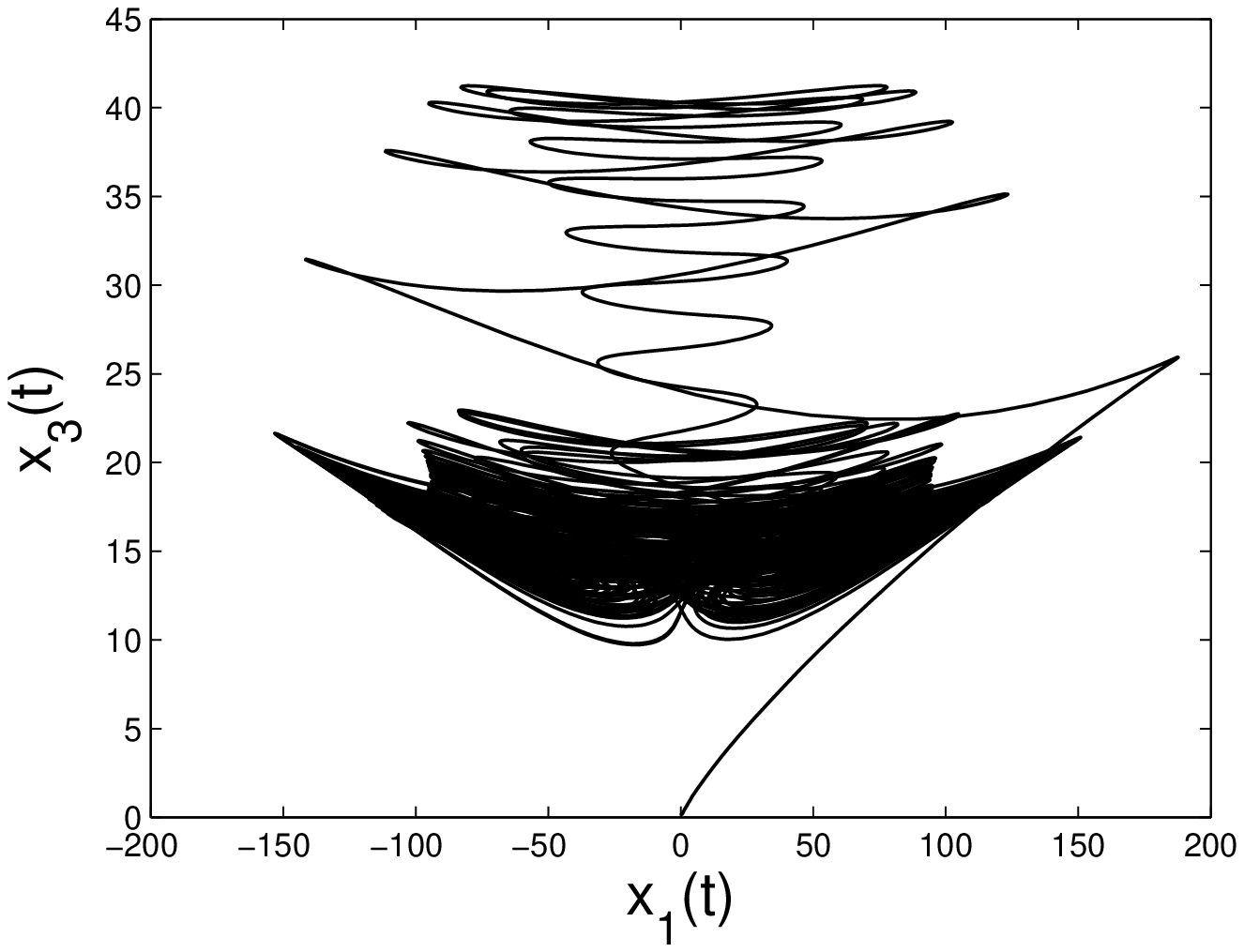}
		\end{minipage}
		\begin{minipage}{200pt}
		\includegraphics[width=200pt]{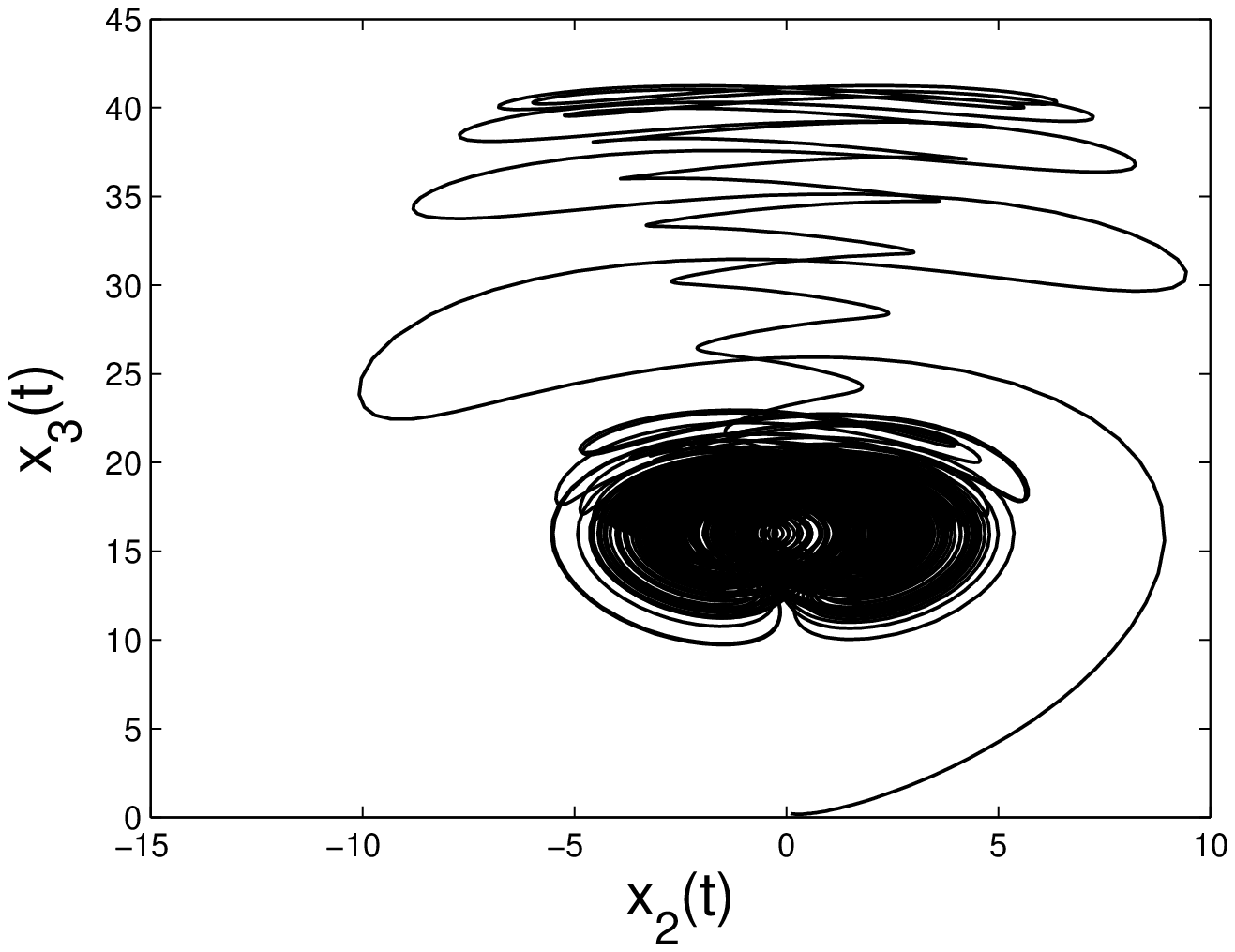}
		\end{minipage}		
	\end{tabular}
    \end{center}
	\caption{Numerical experiment results for chaotic system (\ref{exp1}) obtained at instance of $\alpha=0.95$.}\label{Fig1b}
\end{figure}

\begin{figure}[!h]
	\centering
	\begin{tabular}{cc}
		\begin{minipage}{200pt}
			\includegraphics[width=200pt]{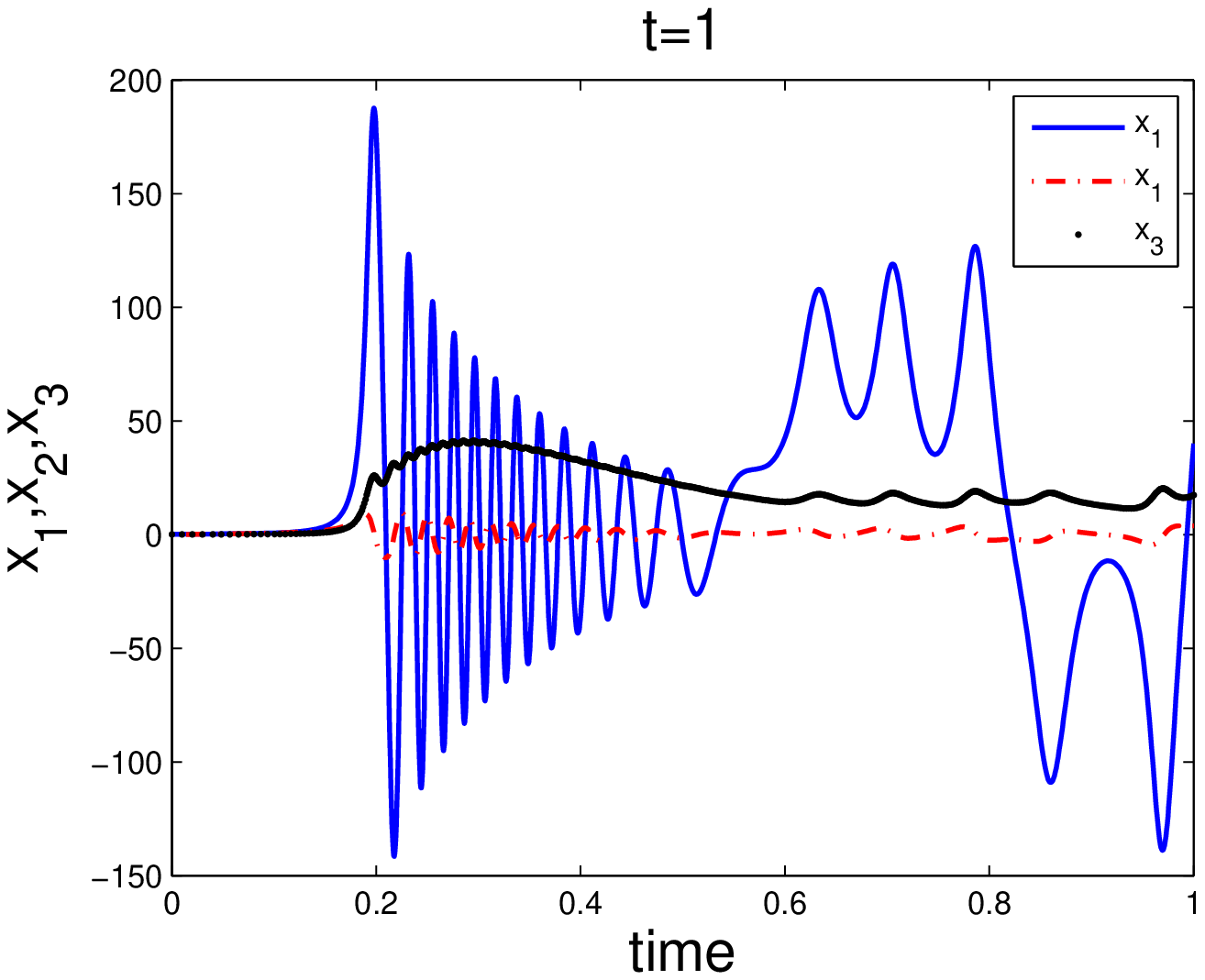}
		\end{minipage}
		\begin{minipage}{200pt}
			\includegraphics[width=200pt]{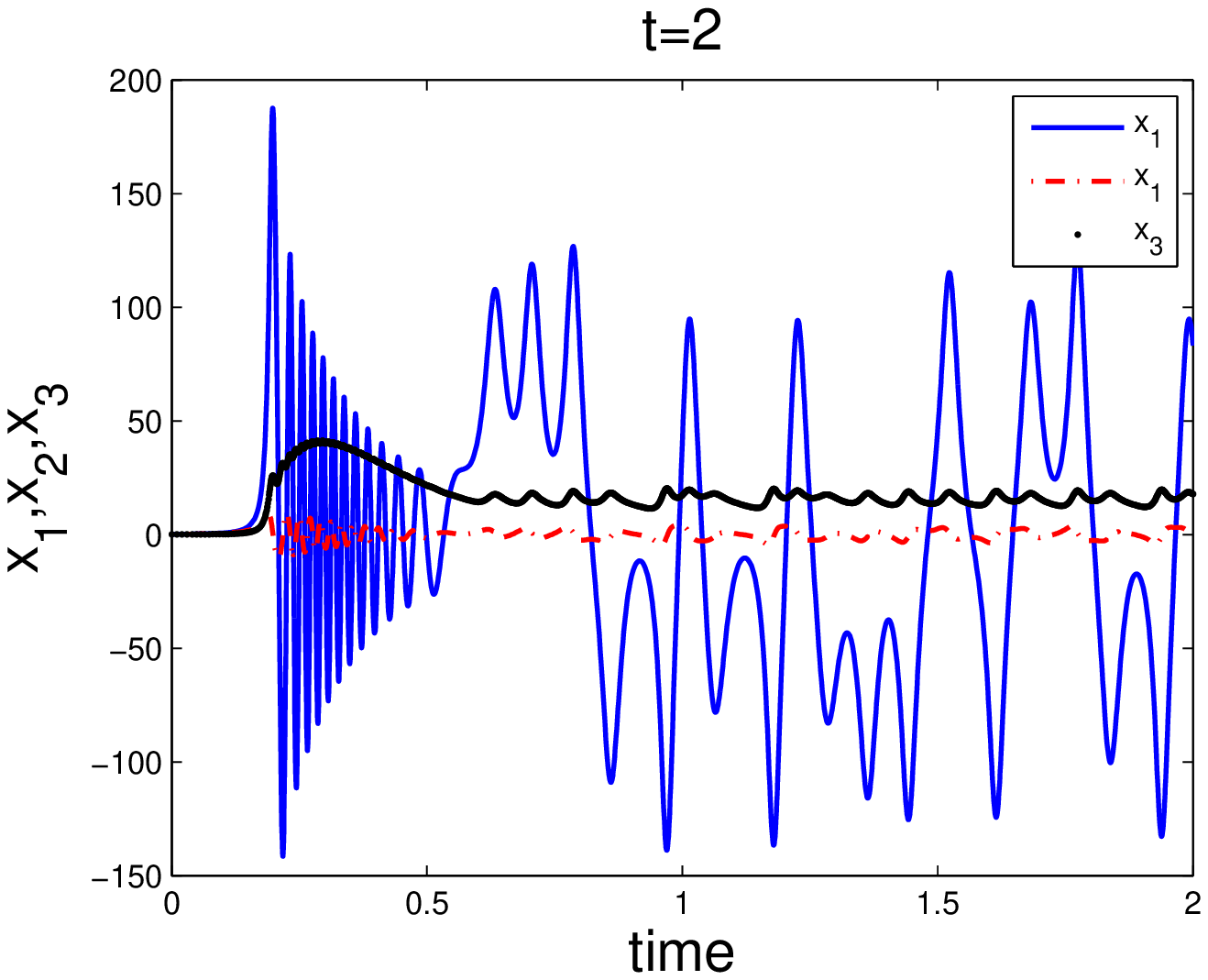}
		\end{minipage}\\
        \begin{minipage}{200pt}
			\includegraphics[width=200pt]{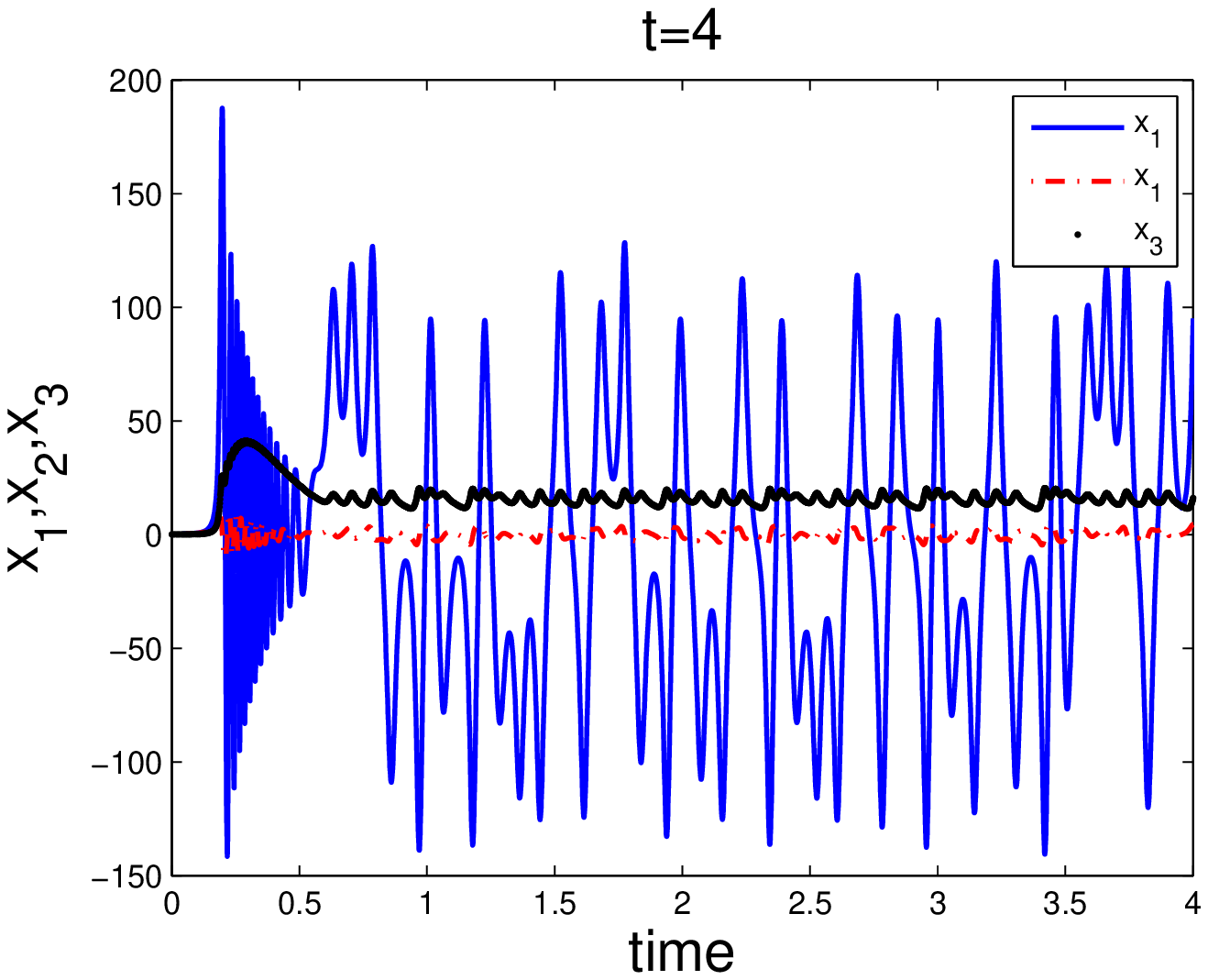}
		\end{minipage}
		\end{tabular}
	\caption{Time-series result for system (\ref{exp1}). The three species interact chaotically in phase as a function of time at $\alpha=0.89$.}\label{Fig1c}
\end{figure}
The novel three dimensional system (\ref{exp1}) shows a highly chaotic and oscillatory attractor for the choice of parameters values $\phi=12, \varphi=16, \psi=5, \sigma=96, \delta=10$. The 3-D phase portrait and 2-D projections with the initial data $x_1(0)=0.2, x_2(0)=0.1, x_3(0)=0.2$ for two instances of fractional index $\alpha=0.75$ and $\alpha=0.95$ as displayed in Figures \ref{Fig1a} and \ref{Fig1b}, respectively. Numerical simulation is performed with time-step $h=0.01$ and final computational time $t=100$. The chaotic and spatiotemporal oscillation of $(x_1,x_2,x_3)$ at different instances of time $(t=1,2,4)$ is shown in Figure \ref{Fig1c}.

\subsection{Example 2} 
\begin{figure}[!h]
	\begin{center}
		\begin{tabular}{cc}
			\begin{minipage}{200pt}
				\includegraphics[width=200pt]{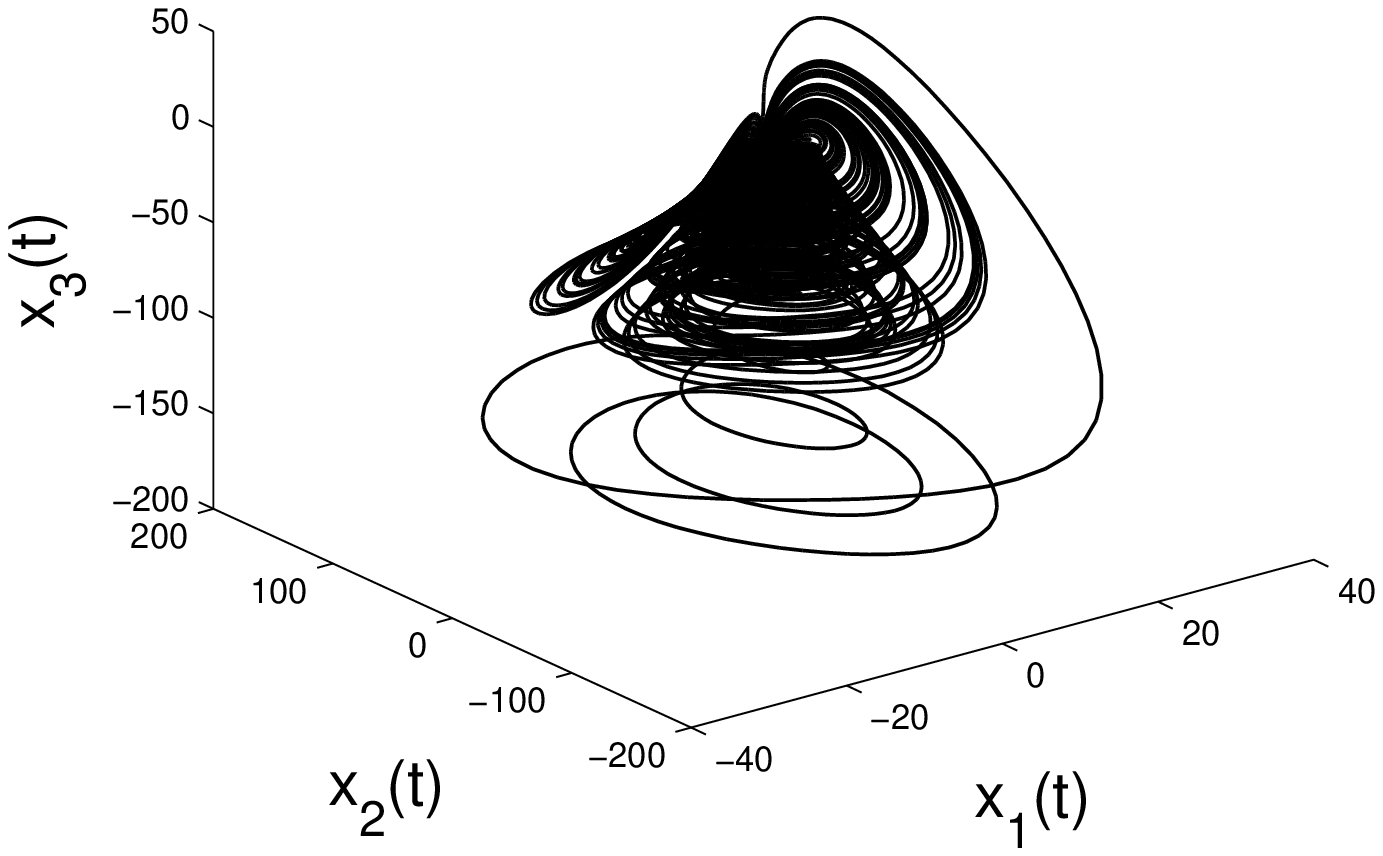}
			\end{minipage}
			\begin{minipage}{200pt}
				\includegraphics[width=200pt]{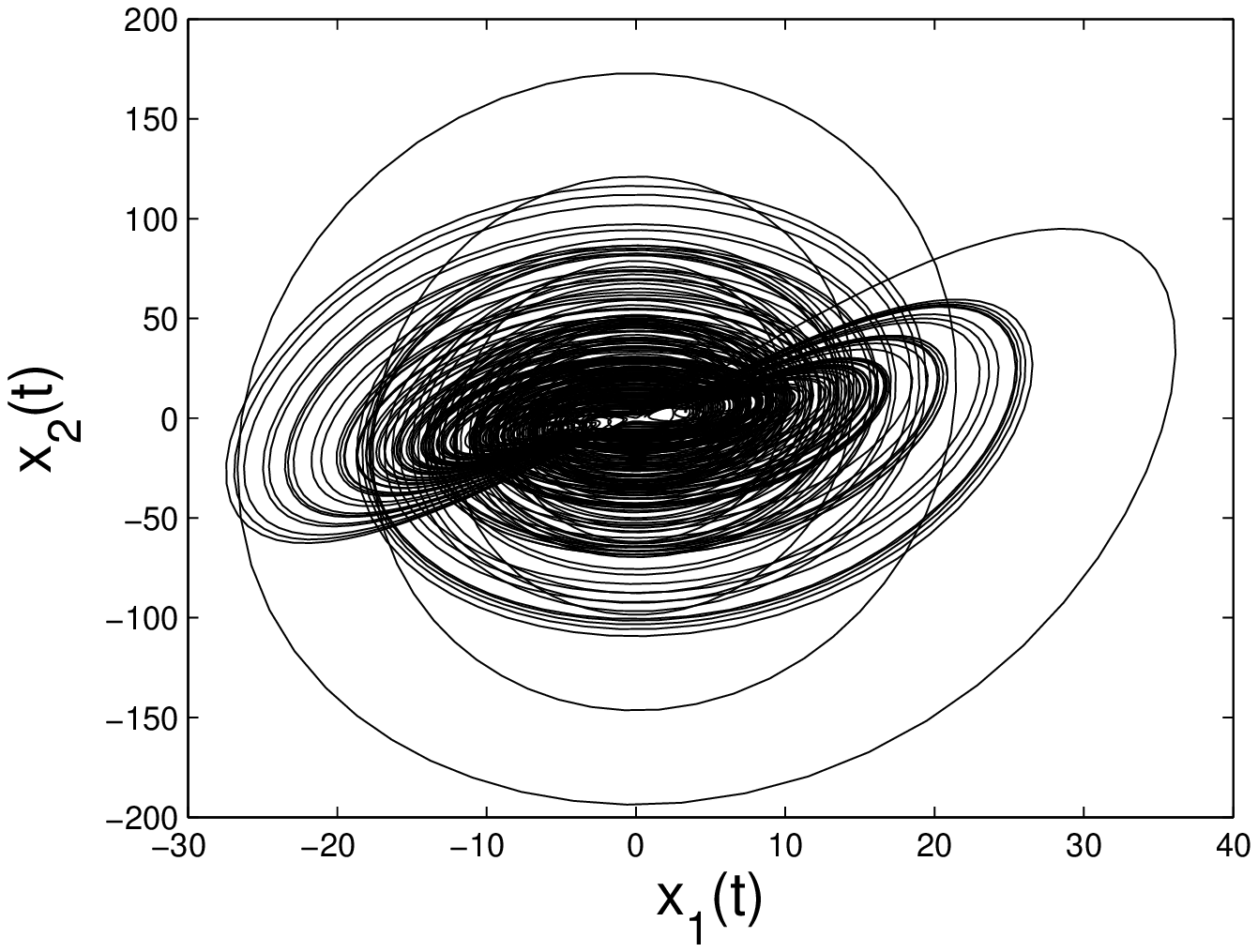}
			\end{minipage}\\
			\begin{minipage}{200pt}
				\includegraphics[width=200pt]{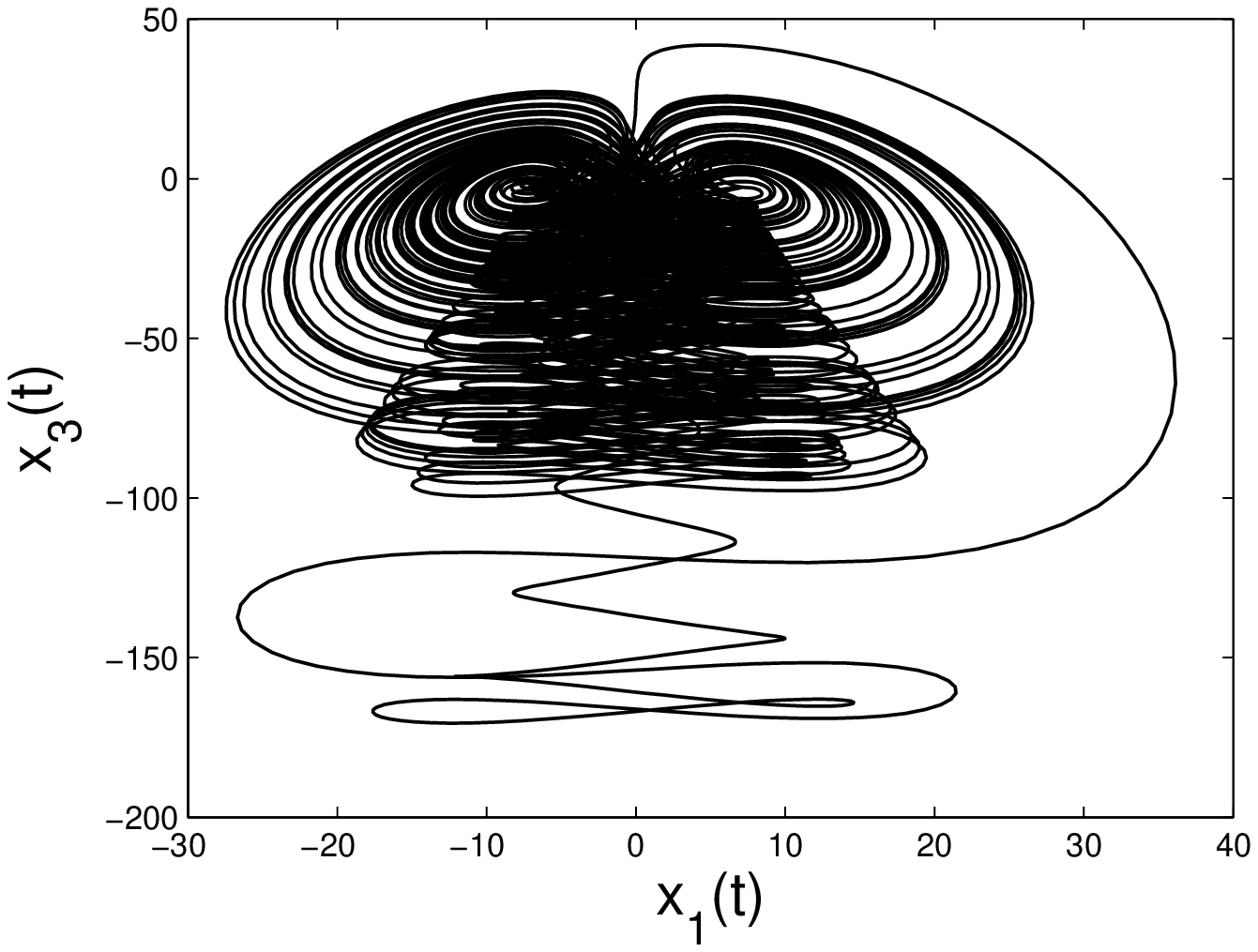}
			\end{minipage}
			\begin{minipage}{200pt}
				\includegraphics[width=200pt]{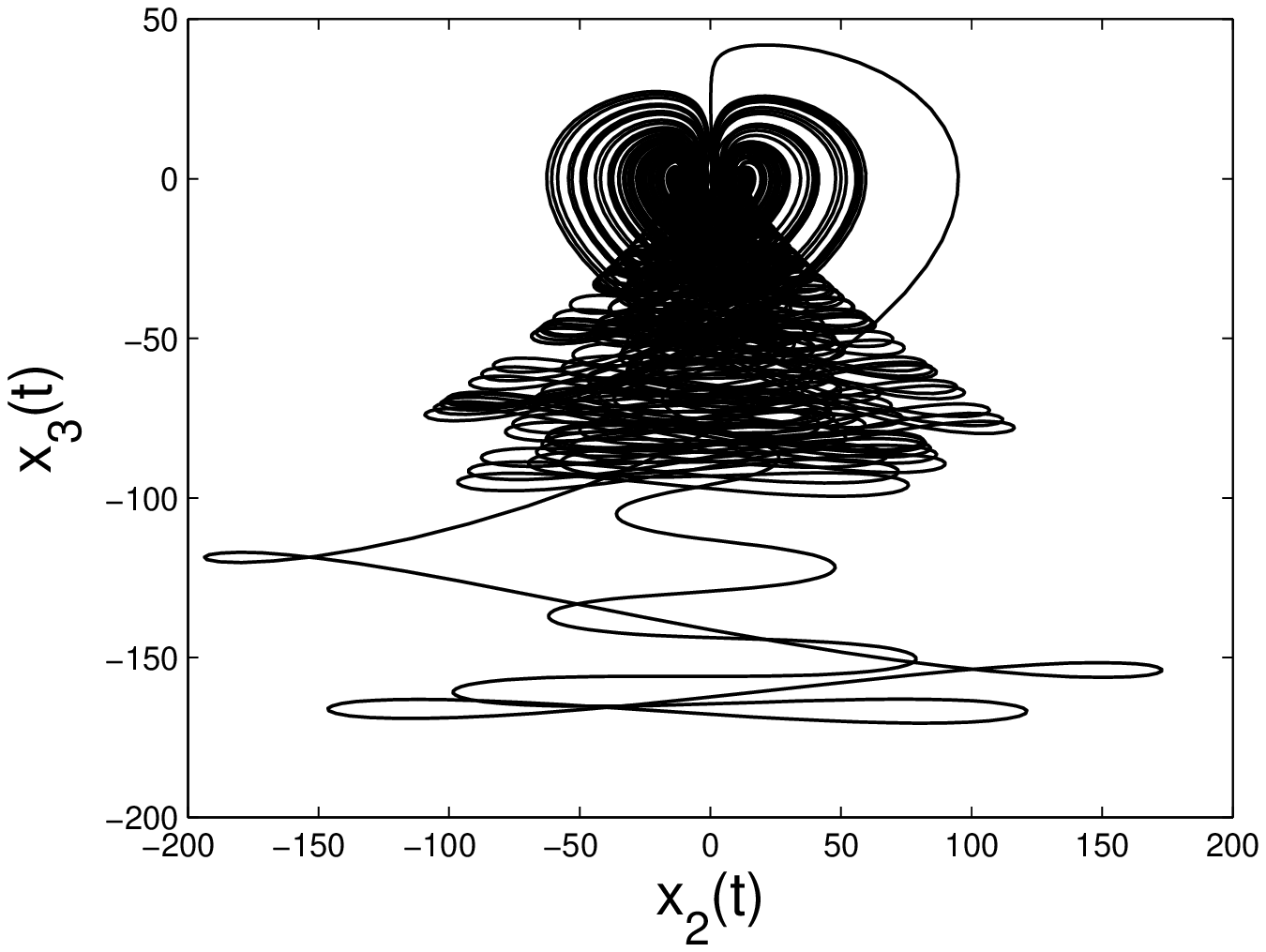}
			\end{minipage}		
		\end{tabular}
	\end{center}
	\caption{Numerical simulation results for chaotic system (\ref{exp2}) obtained at instance of $\alpha=0.75$. Plots  depict the species strange attractor in $\mathbb{R}^3$, and  2-D phase portraits.}\label{Fig2a}
\end{figure}

\begin{figure}[!h]
	\begin{center}
		\begin{tabular}{cc}
			\begin{minipage}{200pt}
				\includegraphics[width=200pt]{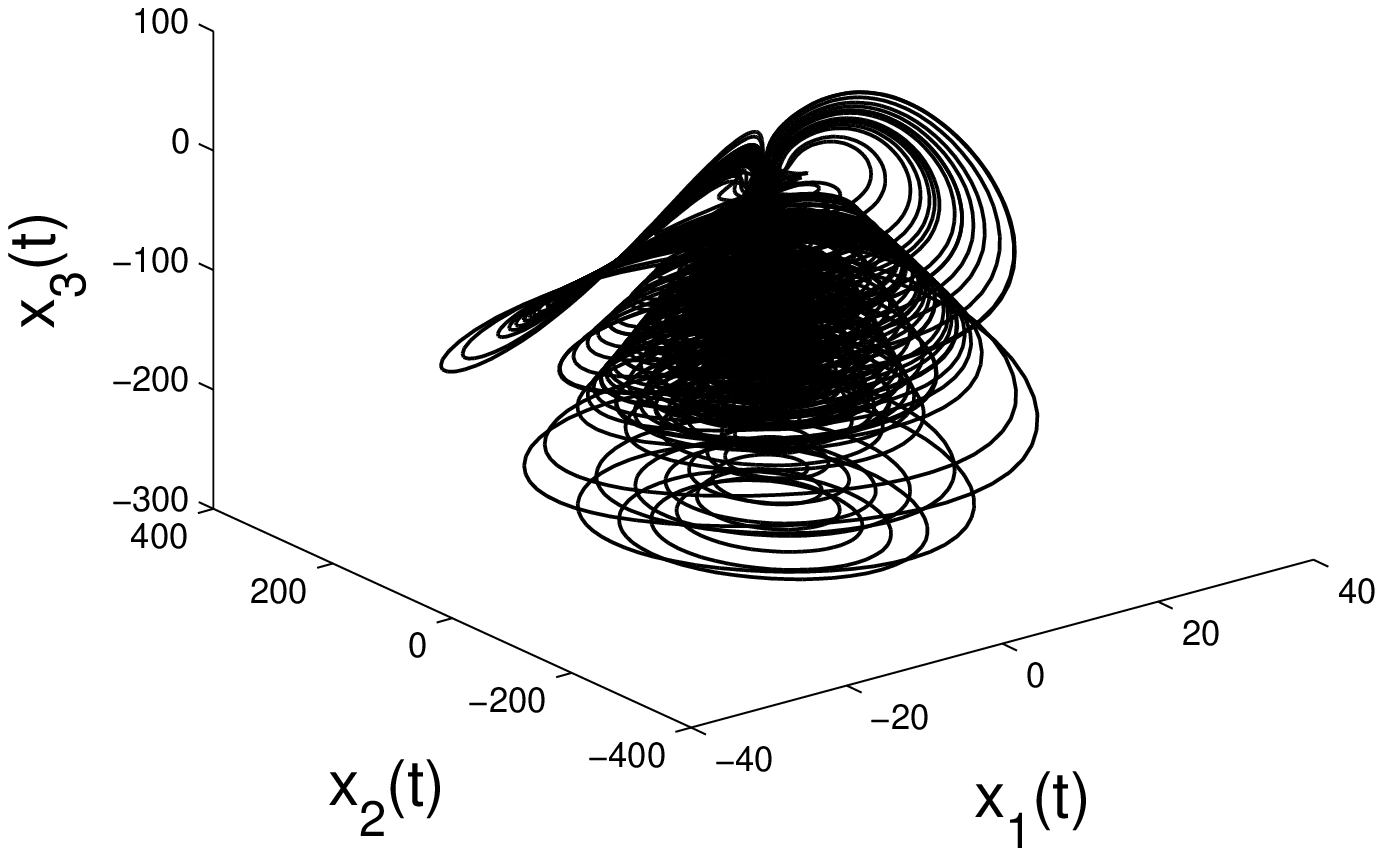}
			\end{minipage}
			\begin{minipage}{200pt}
				\includegraphics[width=200pt]{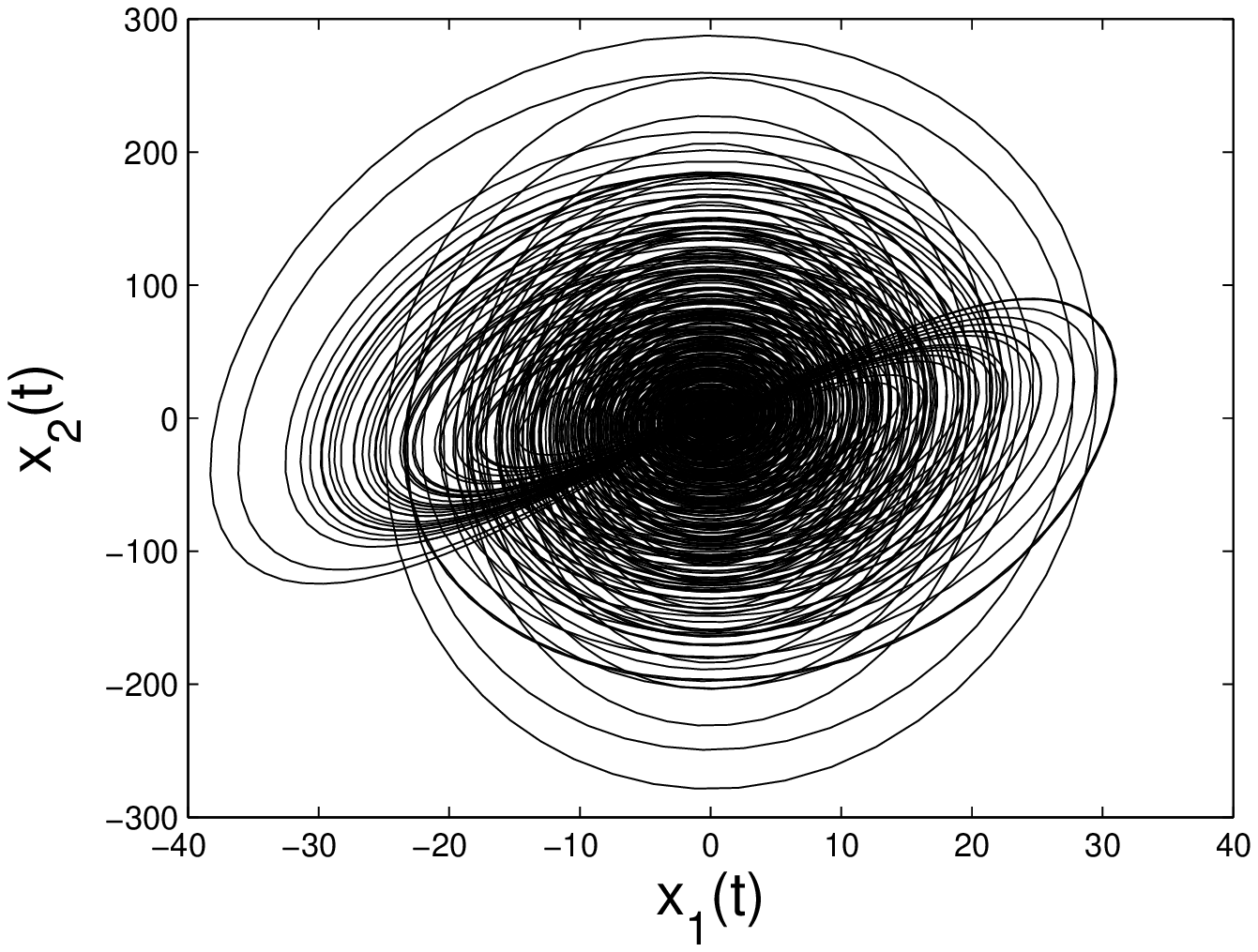}
			\end{minipage}\\
			\begin{minipage}{200pt}
				\includegraphics[width=200pt]{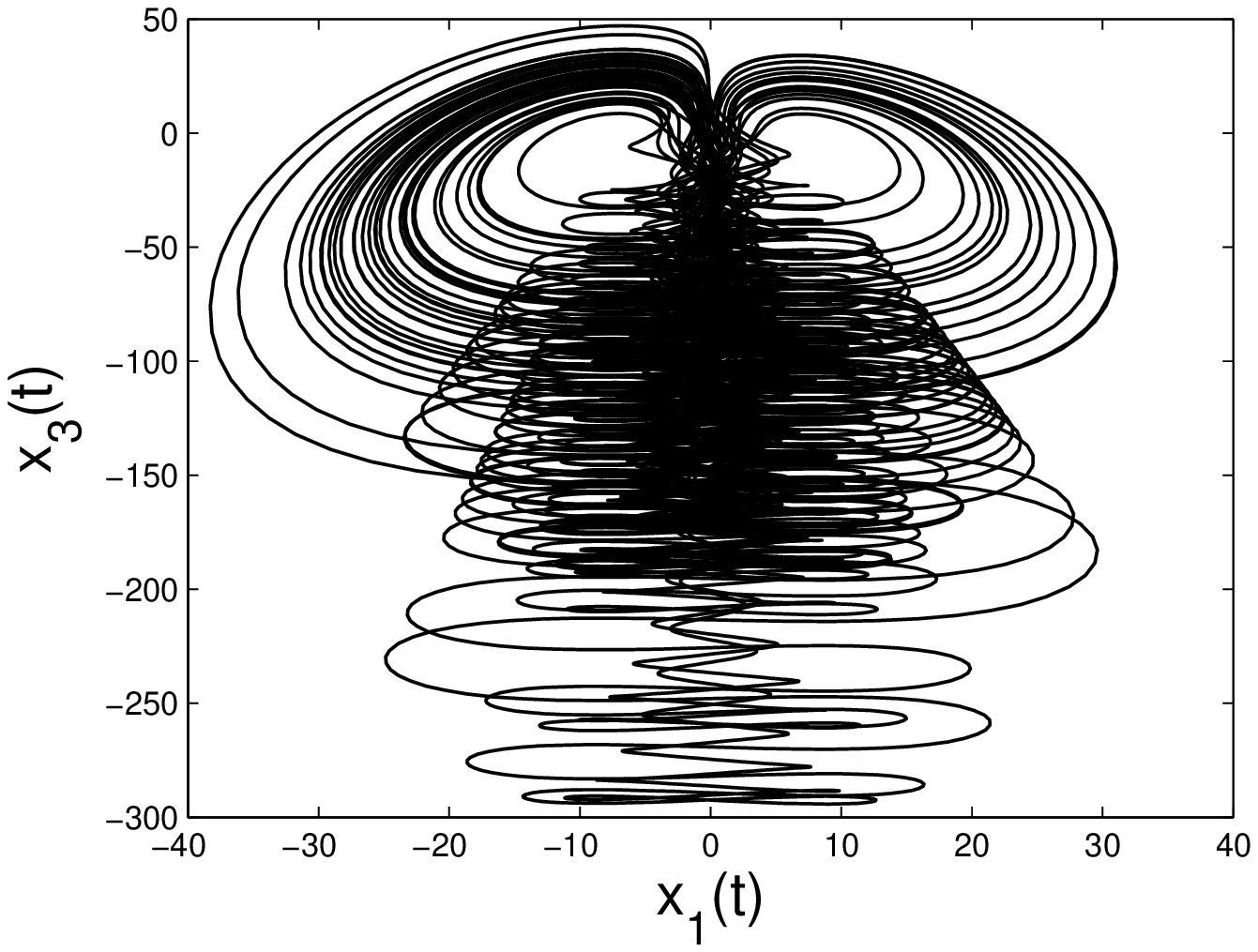}
			\end{minipage}
			\begin{minipage}{200pt}
				\includegraphics[width=200pt]{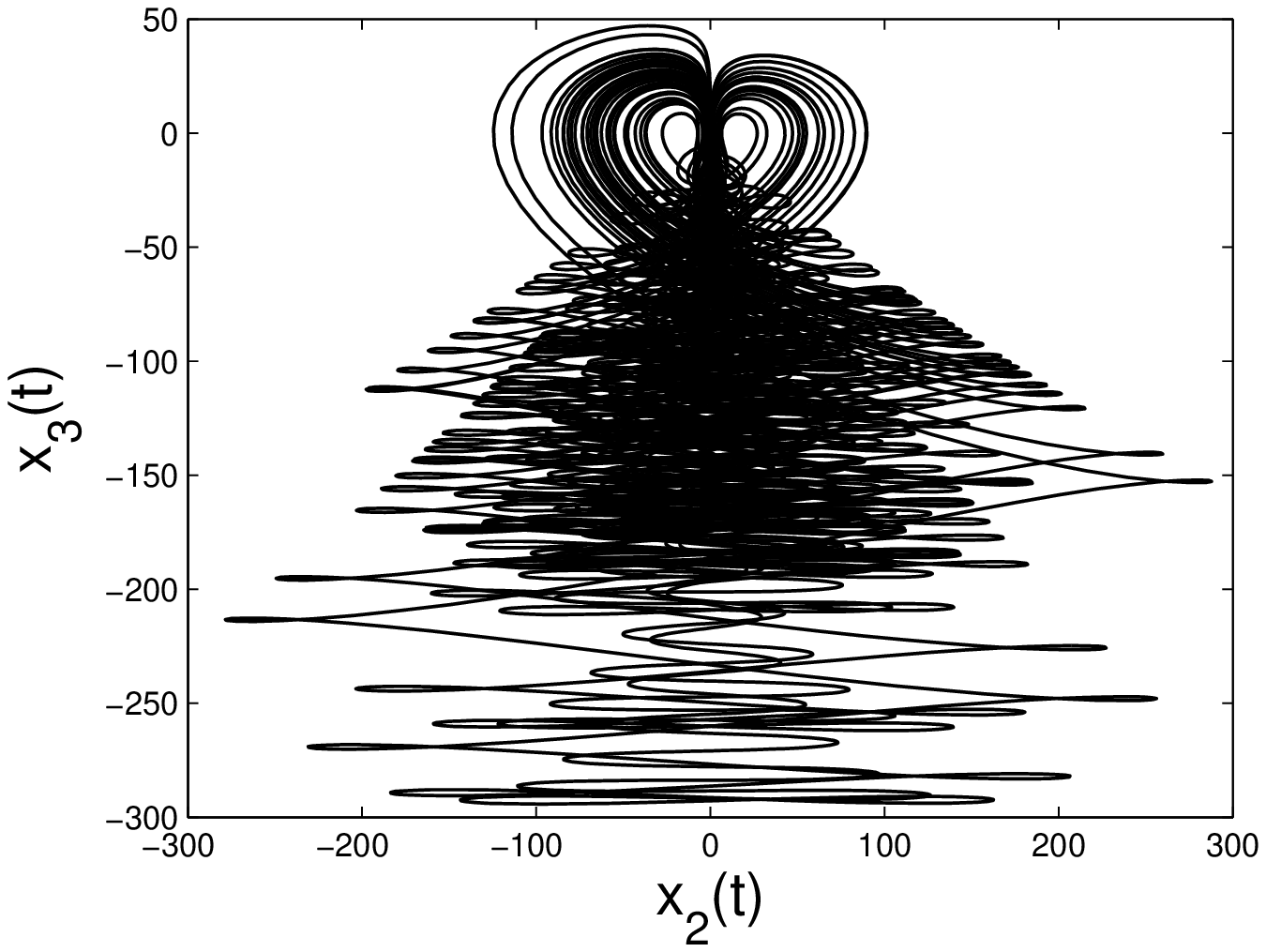}
			\end{minipage}		
		\end{tabular}
	\end{center}
	\caption{Numerical simulation results for chaotic system (\ref{exp2}) obtained at instance of $\alpha=0.95$. Results  represent the species strange attractor in $\mathbb{R}^3$, and  2-D phase portraits.}\label{Fig2b}
\end{figure}

\begin{figure}[!h]
	\begin{center}
	\begin{tabular}{cc}
		\begin{minipage}{200pt}
		\includegraphics[width=200pt]{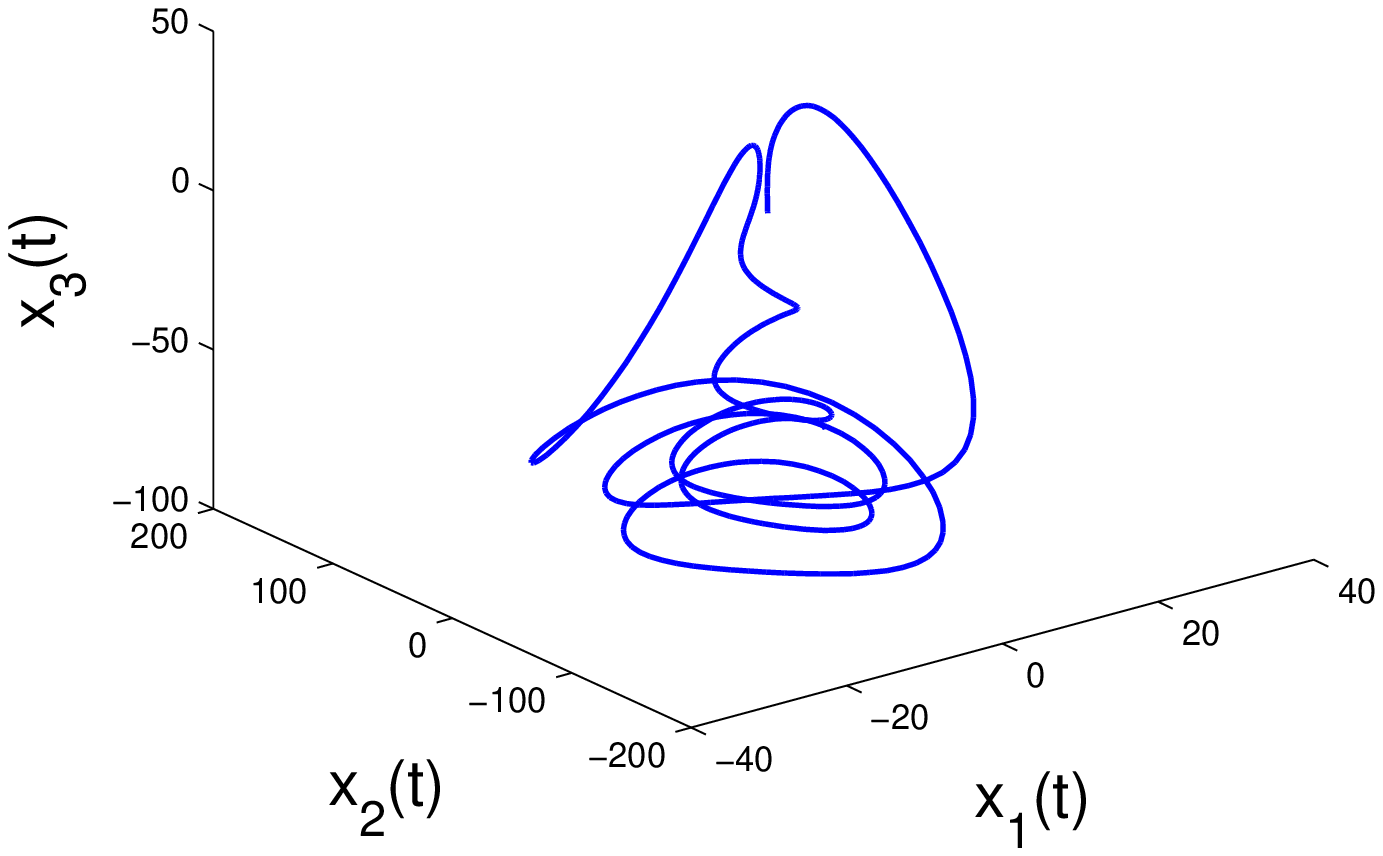}
		\end{minipage}
		\begin{minipage}{200pt}
		\includegraphics[width=200pt]{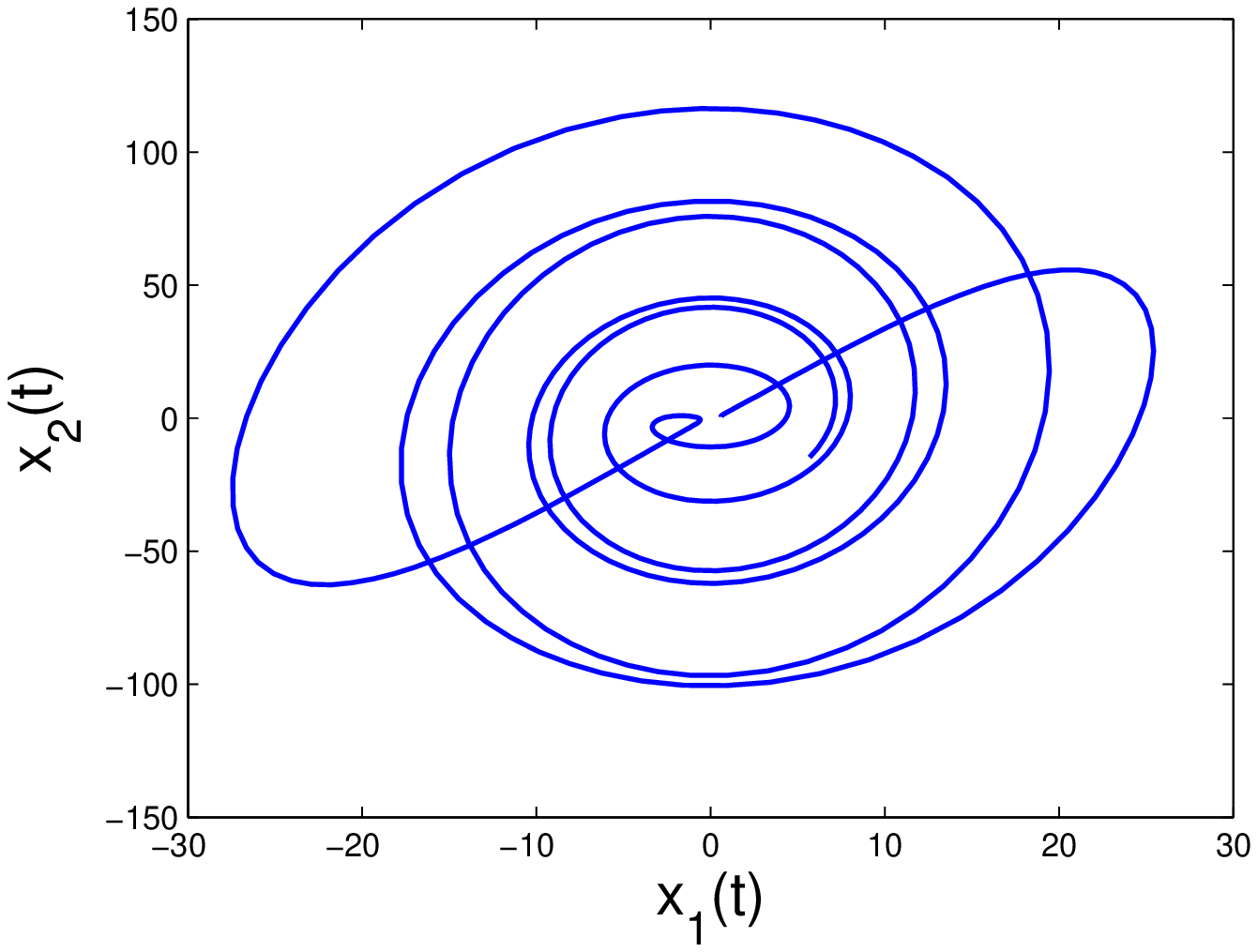}
		\end{minipage}\\
		\begin{minipage}{200pt}
		\includegraphics[width=200pt]{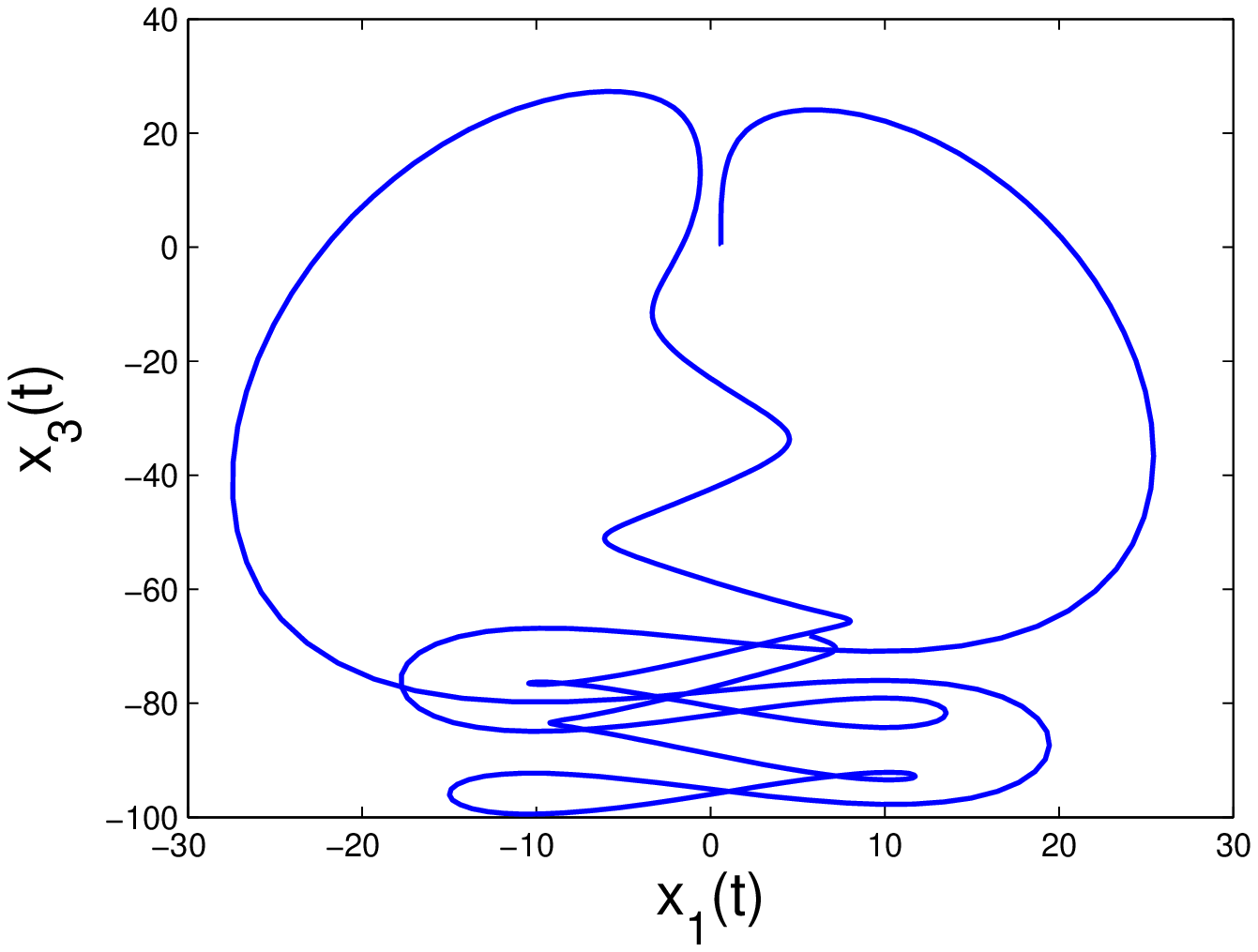}
		\end{minipage}
		\begin{minipage}{200pt}
		\includegraphics[width=200pt]{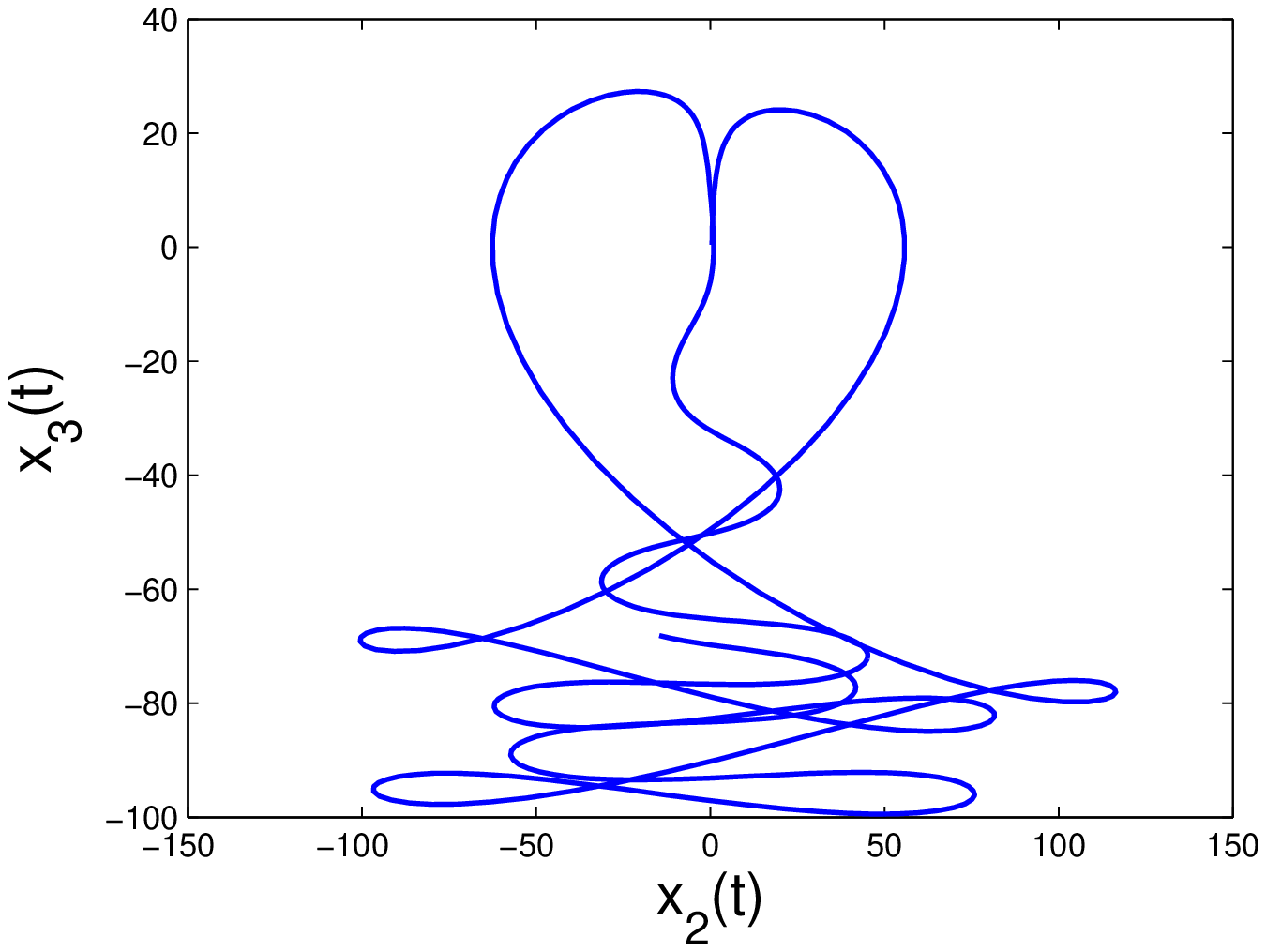}
		\end{minipage}			
	\end{tabular}
    \end{center}
	\caption{Emergence of some strange chaotic attractors for system (\ref{exp2}) at $t=5 $ and $\alpha=0.68$. Other parameters are as given above.}\label{Fig2c}
\end{figure}

For the second example, consideration is given to a six-term novel 3-D dissipative chaotic system \cite{Aza16} which is illustrated by the fractional 3-D dynamics
\begin{equation}\label{exp2}
\begin{split}
^{ABC}_0\mathcal{D}_t^\alpha x_1(t)=&f_1(x_1,x_2,x_3)=\phi(x_2(t)-x_1(t)),\\
^{ABC}_0\mathcal{D}_t^\alpha x_2(t)=&f_2(x_1,x_2,x_3)= x_1(t)-x_1(t) x_3(t),\\
^{ABC}_0\mathcal{D}_t^\alpha x_3(t)=&f_3(x_1,x_2,x_3)=50-\varphi x_1^2(t)-\psi x_3(t),
\end{split}
\end{equation}
where $x_1(t), x_2(t), x_3(t)$ are the densities and $\phi>0,\varphi>0,\psi>0$ are parameters.
The simulation results in Figures \ref{Fig2a} and \ref{Fig2b} are obtained with $\alpha=0.75$ and $\alpha=0.95$, respectively. We utilize parameter values $\phi=2.6, \varphi=0.5, \psi=0.4$, subject to the initial conditions $x_1(0)=0.6, x_2(0)=0.5$ and $x_3(0)=0.4$. The Initial formation of chaotic and spatiotemporal patterns obtained in Figure \ref{Fig2c} corresponds to $\alpha=0.68$ and time $t=5$. It should be noted that the species oscillate chaotically in phase regardless of the variation in $\alpha$ and computational time $t$.

\begin{figure}[!ht]
	\centering
	\begin{tabular}{cc}
		\begin{minipage}{200pt}
			\includegraphics[width=200pt]{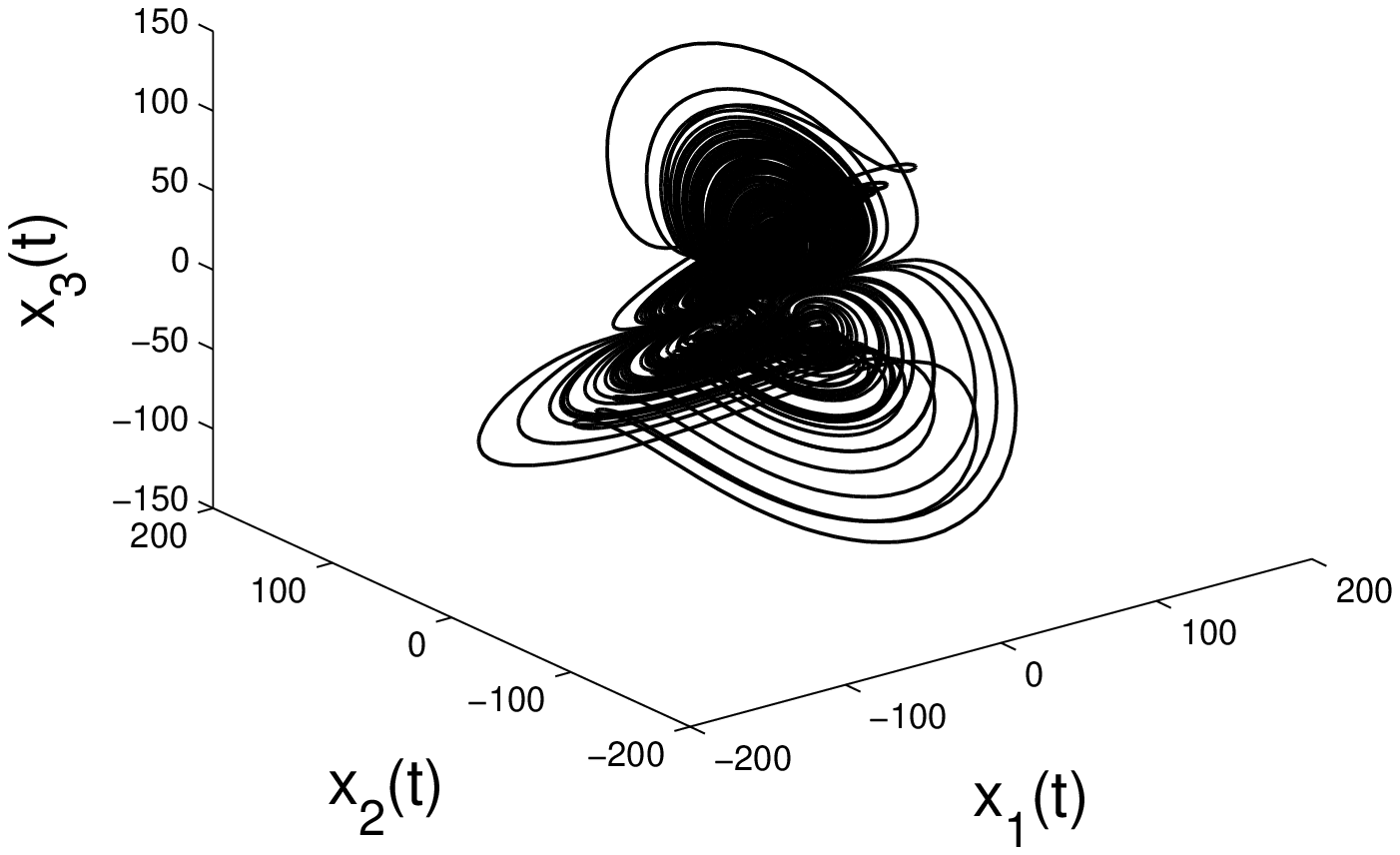}
		\end{minipage}
		\begin{minipage}{200pt}
			\includegraphics[width=200pt]{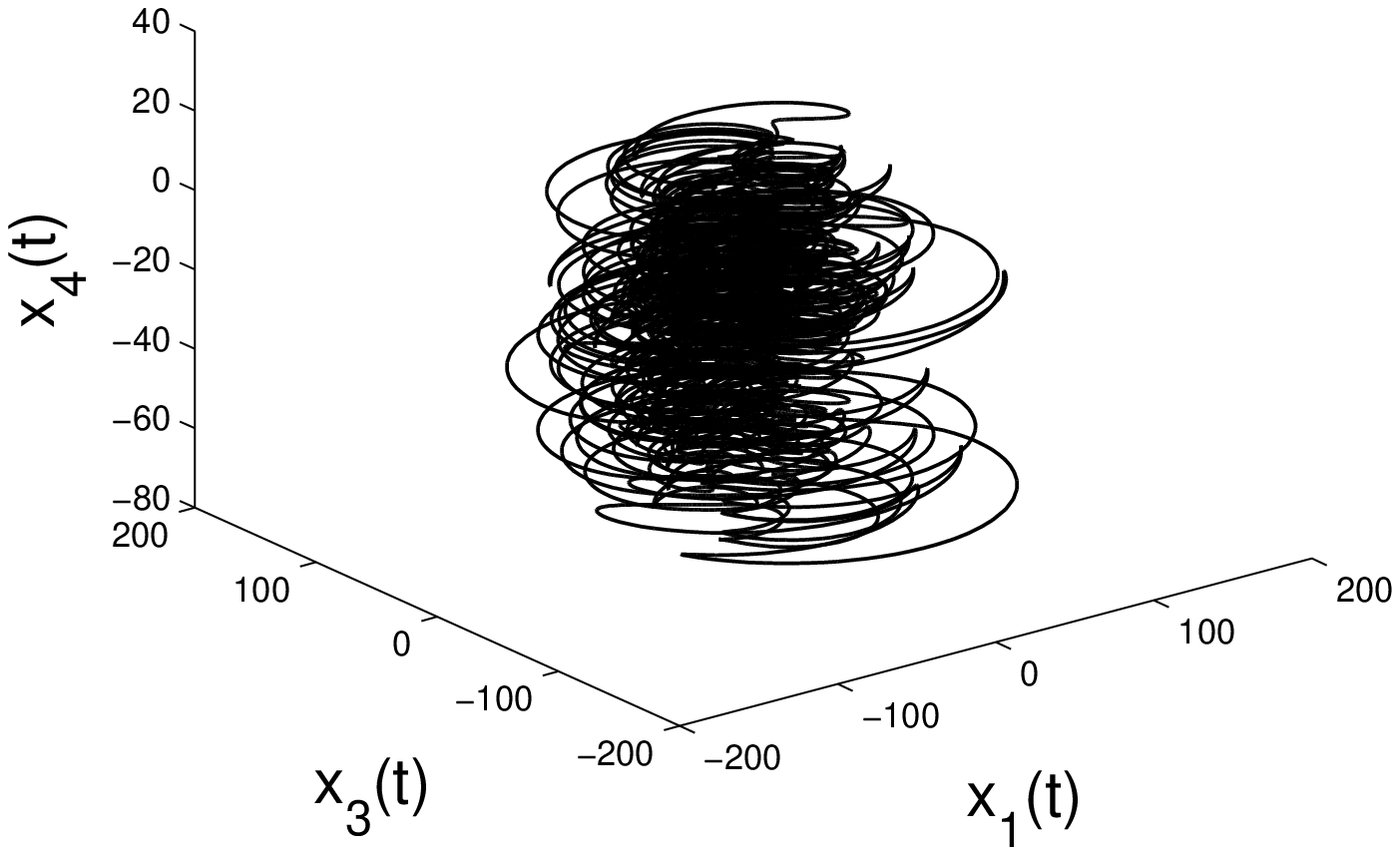}
		\end{minipage}\\
  		\begin{minipage}{200pt}
			\includegraphics[width=200pt]{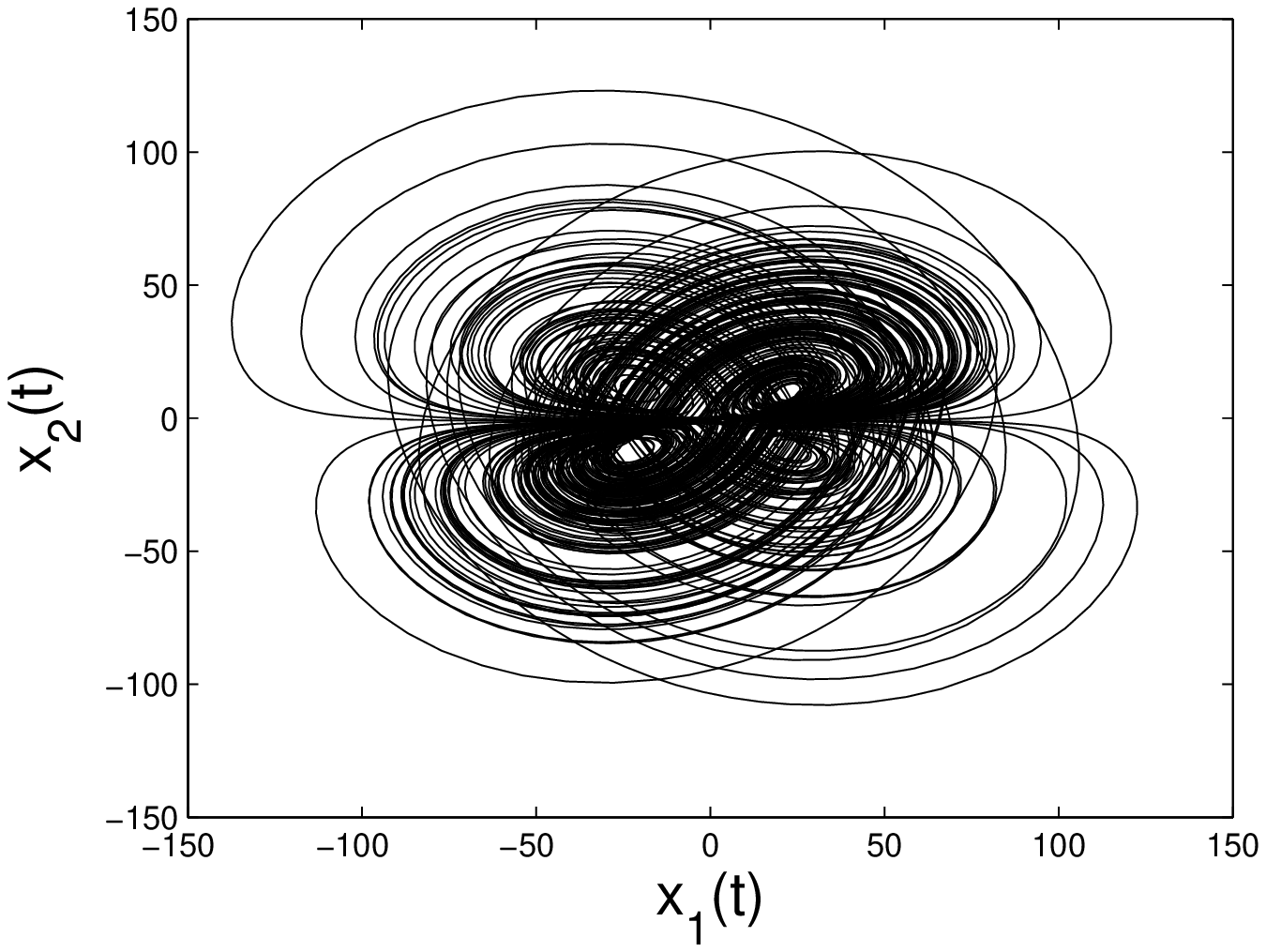}
		\end{minipage}
		\begin{minipage}{200pt}
			\includegraphics[width=200pt]{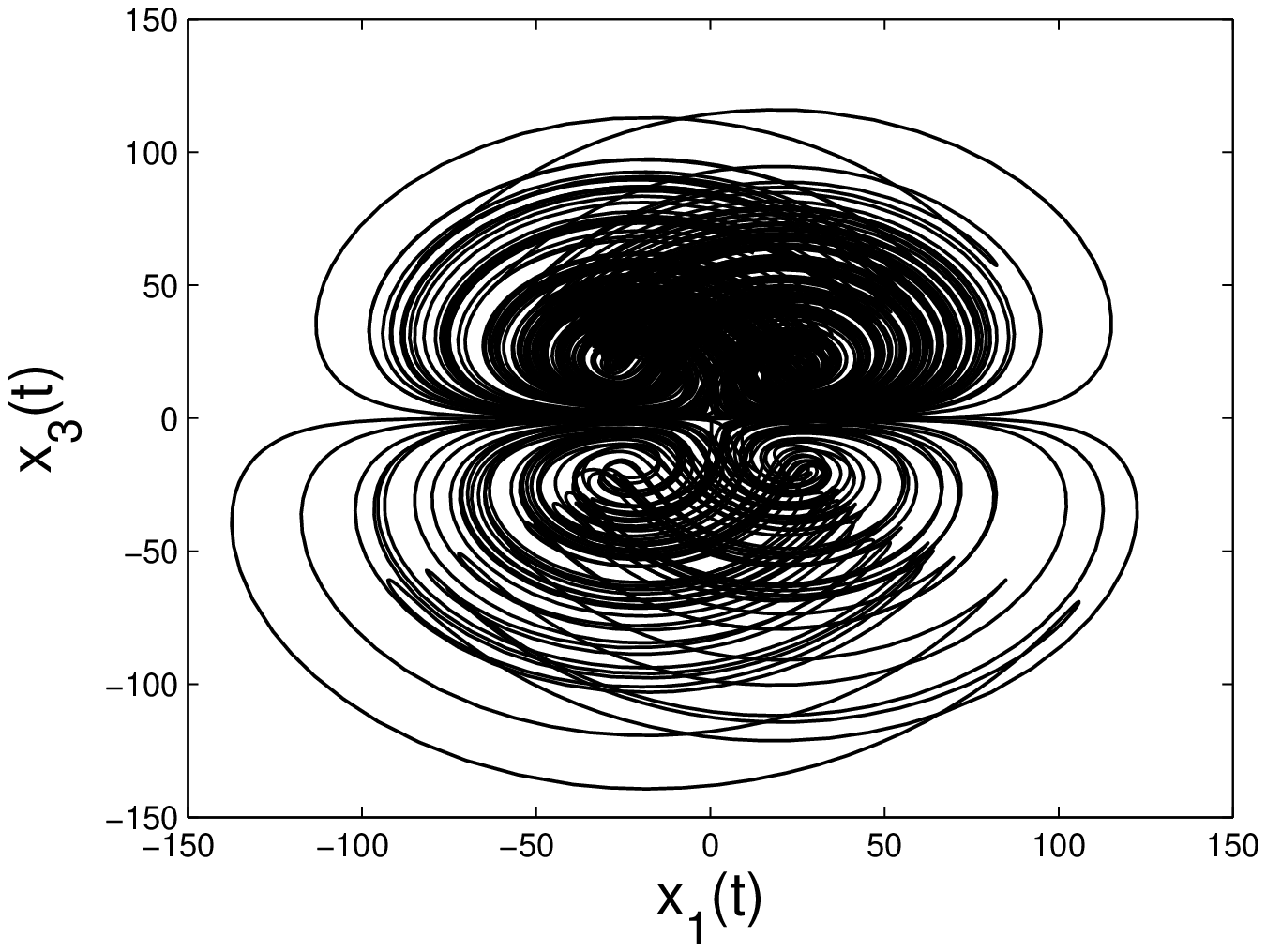}
		\end{minipage}\\
      		\begin{minipage}{200pt}
			\includegraphics[width=200pt]{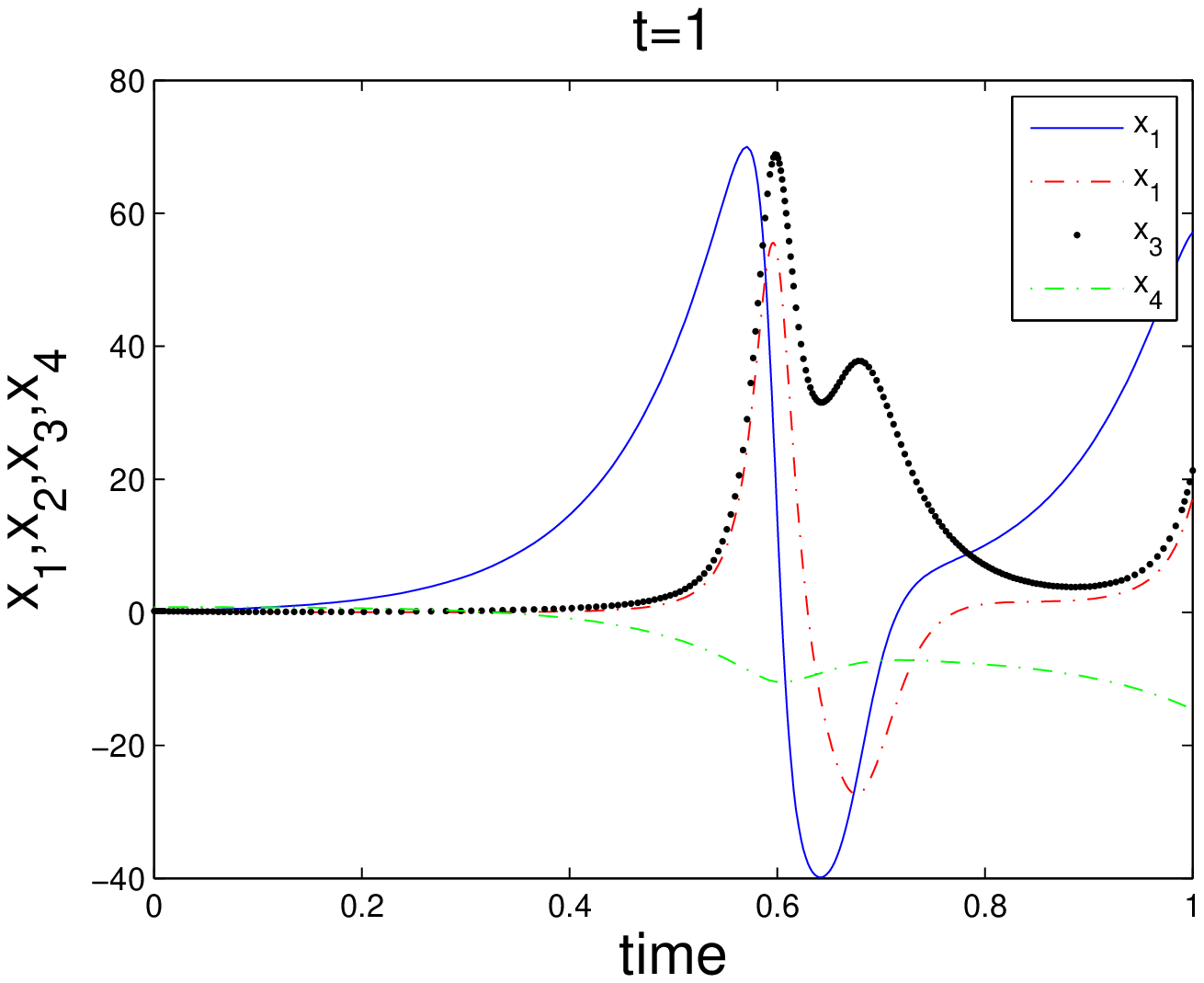}
		\end{minipage}
		\begin{minipage}{200pt}
			\includegraphics[width=200pt]{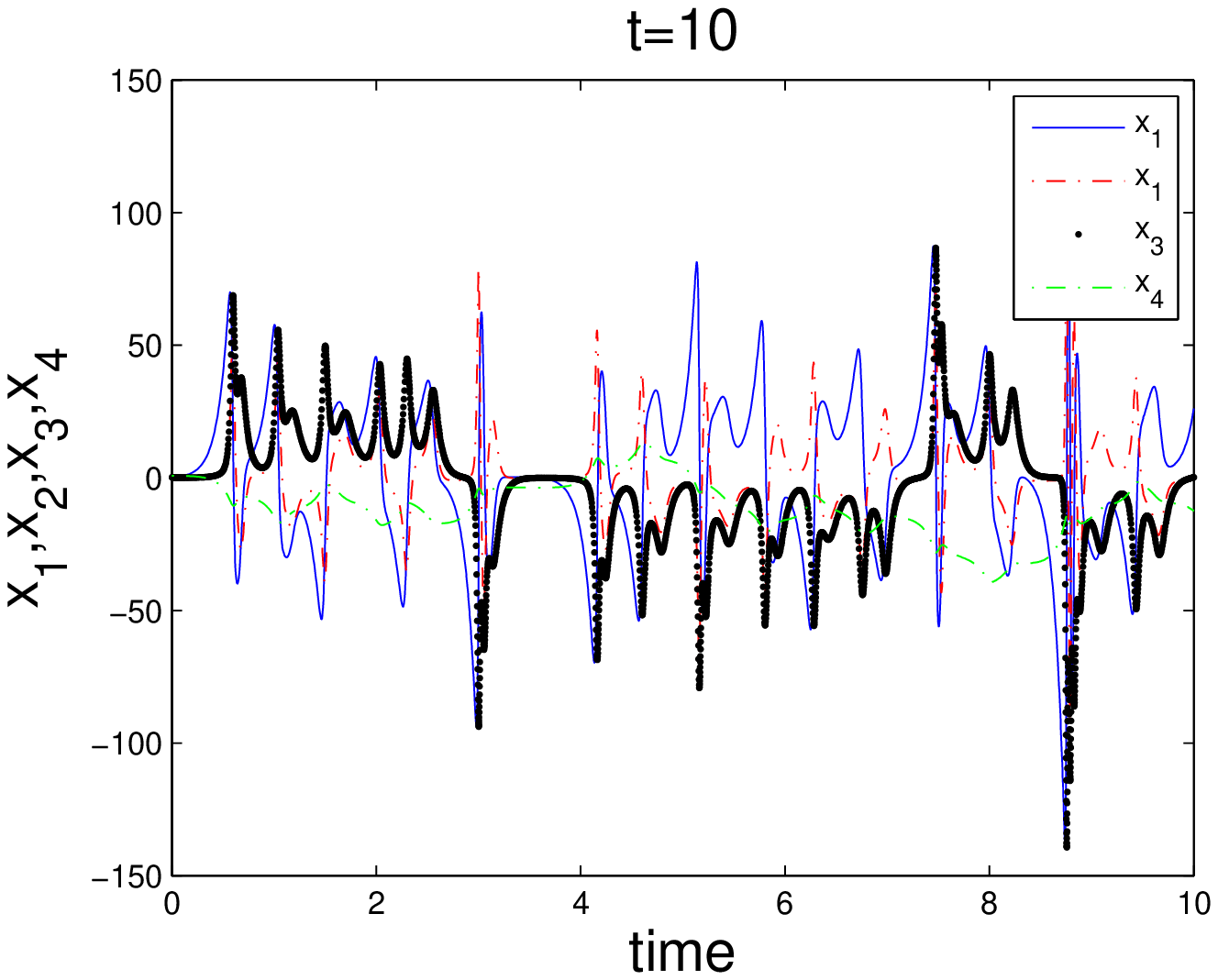}
		\end{minipage}
		\end{tabular}
	\caption{Numerical simulation of fractional system (\ref{exp3}) showing hyperchaotic oscillations in 3-D $(upper-row)$ and 2-D $(lower-row)$ projections in phase. The lower-row corresponds to time series solutions at $t=(1,10)$ for $\alpha=0.93$.}\label{Fig3}
\end{figure}

\subsection{Example 3} We extend our numerical experiment to a four-component system, by considering the novel 4-D hyper chaotic four-wing fractional equations
 \begin{equation}\label{exp3}
\begin{split}
^{ABC}_0\mathcal{D}_t^\alpha x_1(t)=&f_1(x_1,x_2,x_3,x_4)=\phi x_1(t)-x_2(t)x_3(t)+x_4(t),\\
^{ABC}_0\mathcal{D}_t^\alpha x_2(t)=&f_2(x_1,x_2,x_3,x_4)= -\varphi x_2(t)+ x_1(t) x_3(t)+x_4(t),\\
^{ABC}_0\mathcal{D}_t^\alpha x_3(t)=&f_3(x_1,x_2,x_3,x_4)=-\psi x_3(t)+ x-1(t)x_1^2(t)+ x_1(t),\\
^{ABC}_0\mathcal{D}_t^\alpha x_4(t)=&f_1(x_1,x_2,x_3,x_4)=\sigma x_1(t),
\end{split}
\end{equation}
where $x_1(t), x_2(t), x_3(t)$ and $x_4(t)$ are the densities and $\phi>0,\varphi>0,\psi>0, \sigma>0$ are mere constant parameters.
By applying (\ref{comp}), equation (\ref{exp2}) can be written in the form
\begin{eqnarray*}
x_{1,n+1}&=&x_{1,n}+\omega_1(n,\alpha,h)f_1(t_n,x_{1,n},x_{2,n},x_{3,n},x_{4,n}) +\omega_2(n,\alpha,h)\\
&&\times f_1(t_{n-1},x_{1,n-1},x_{2,n-1},x_{3,n-1},x_{4,n-1})\nonumber\\
x_{2,n+1}&=&x_{2,n}+\omega_1(n,\alpha,h)f_2(t_n,x_{1,n},x_{2,n},x_{3,n},x_{4,n}) +\omega_2(n,\alpha,h)\\
&&\times f_2(t_{n-1},x_{1,n-1},x_{2,n-1},x_{3,n-1},x_{4,n-1})\nonumber\\
x_{3,n+1}&=&x_{3,n}+\omega_1(n,\alpha,h)f_3(t_n,x_{1,n},x_{2,n},x_{3,n},x_{4,n}) +\omega_2(n,\alpha,h)\\
&&\times f_3(t_{n-1},x_{1,n-1},x_{2,n-1},x_{3,n-1},x_{4,n-1})\nonumber\\
x_{4,n+1}&=&x_{4,n}+\omega_1(n,\alpha,h)f_4(t_n,x_{1,n},x_{2,n},x_{3,n},x_{4,n}) +\omega_2(n,\alpha,h)\\
&&\times f_4(t_{n-1},x_{1,n-1},x_{2,n-1},x_{3,n-1},,x_{4,n-1})\nonumber
\end{eqnarray*}

 The 4-D chaotic results given in Figure \ref{Fig3} is obtained with the parameter values
$$\phi=8,\;\;\varphi=33,\;\;\psi=16,\;\; \sigma=1.25.$$
We use the initial data $x_1(0)=0.2, x_2(0)=0.4, x_3(0)=0.2$  and $x_4(0)=0.7$ for the numerical simulation at $\alpha=0.93$ with final computational time $t=150$ for the upper- and middle-rows. The 3-D hyperchaotic projections is obtained on $(x_1(t),x_2(t),x_3(t)), (x_1(t),x_3(t),x_4(t))$ and $(x_2(t),x_3(t),x_4(t))$ respectively. The 3-D  and 2-D projections have the shape of four-wind attractors. The result displayed in the bottom-row depicts time series solution, obtained at instances of time $t=(1,10)$.

\section{Conclusion}
Mathematical analysis and simulation of the newly proposed Adams-Bashforth scheme with the Atangana-Baleanu in the sense of Caputo, when applied to study different chaotic and hyperchaotic systems that are of current and recurrent interests are considered in this paper. We considered a general two-component fractional differential system to study the existence and uniqueness of solutions. Numerical results which reveal some strange attractors at some instances of fractional power $\alpha$ are given to address the points and queries therein. In the future work, the existence and uniqueness of solution reported for a general two-component differential equations will be extended to multi-dimensional problems.


\section*{Correction to: Chaos {\color{blue}https://doi.org/10.1063/1.5085490}.}	

\section*{Introduction}
Following the publication of an article in \cite{Owo19}, we provide corrections to mistakes therein as pointed out by some anonymous reviewers. Details of the original article can be found online at {\color{blue}https://doi.org/10.1063/1.5085490}.

{\bf In abstract}\\
Lines 2-3 $\rightarrow$ We adopt the Banach fixed point theorem...

{\bf Correction to equations (10) and (11)}.\\
Equations (10) and (11) in \cite{Owo19} are corrected as follows.
\begin{equation}
\begin{split}
A_{\alpha,1}&=\frac{\alpha}{\Gamma(\alpha)ABC(\alpha)}\int_{0}^{t_{n+1}}(t_{n+1}-t)^{\alpha-1}\left\{\frac{t-t_{n-1}}{h}f(t_n,u_n) - \frac{t-t_{n}}{h}f(t_{n-1},u_{n-1})\right\}d\tau\\
&=\frac{f(t_n,u_n)}{\Gamma(\alpha)ABC(\alpha)h}\left\{\frac{t^{\alpha+1}_{n+1}}{\alpha+1}-t_{n-1}t^\alpha_{n+1} \right\} -
\frac{f(t_{n-1},u_{n-1})}{h\Gamma(\alpha)ABC(\alpha)}\left\{\frac{t^{\alpha+1}_{n+1}}{\alpha+1}-t_{n}t^\alpha_{n+1} \right\}.
\end{split}
\end{equation}
Also,
\begin{equation}
\begin{split}
A_{\alpha,2}&=\frac{f(t_n,u_n)}{h\Gamma(\alpha)ABC(\alpha)}\left\{\frac{t^{\alpha+1}_{n}}{\alpha+1}-t_{n-1}t^\alpha_{n} \right\} -
\frac{f(t_{n-1},u_{n-1})}{h\Gamma(\alpha)ABC(\alpha)}\left\{\frac{t^{\alpha+1}_{n}}{\alpha+1}-t_{n}^{\alpha+1} \right\}.
\end{split}
\end{equation}
	{\bf \color{red}
	Equation (12) is corrected as
	\begin{equation*}
		\begin{split}
			&u(t_{n+1})-u(t_n)=\frac{1-\alpha}{AB(\alpha)}\left\{ f(t_n,u_n) -f(t_{n-1},u_{n-1}) \right\}	+\frac{f(t_n,u_n)}{h\Gamma(\alpha)AB(\alpha)}\left\{ t^{\alpha+1}_{n+1} - t_nt^\alpha_{n+1} \right\}\\
			&- \frac{f(t_{n-1},u_{n-1})}{h\Gamma(\alpha)AB(\alpha)}\left\{ \frac{t_{n+1}^{\alpha+1}}{\alpha+1}-t_nt_{n+1}^\alpha\right\}
			-\frac{f(t_n,u_n)}{h\Gamma(\alpha)AB(\alpha)}\left\{\frac{t_{n}^{\alpha+1}}{\alpha+1}-t_{n-1}t_{n}^\alpha \right\}\\ &+\frac{f(t_{n-1},u_{n-1})}{h\Gamma(\alpha)AB(\alpha)}\left\{ \frac{t^{\alpha+1}_n}{\alpha+1} -t_n^{\alpha+1} \right\}
		\end{split}	
	\end{equation*}
	where $n \ge 1$ and the $u(0)=u_0$ is given by initial condition, where $u_1$ can be obtained using the Runge-Kutta method. }

Under equation (15) in \cite{Owo19}, we assume that the function $f(x,y,t)$ and $g(x,y,t)$ are uniformly Lipschitz continuous in $x$ and $y$ respectively, also $g$ and $f$ are continuous in $t$.

To avoid confusion on $\Gamma$ as gamma function, and $\Gamma$ as defined mapping, we suggest that $\Gamma$ should be replaces by $\Pi$
$$\Pi:G_{a,b}\longrightarrow G_{a,b}$$
$$\Pi x(t)=x_0+\frac{1-\alpha}{AB(\alpha)}f(x,y,t)+\frac{\alpha}{AB(\alpha)\Gamma(\alpha)}\int_{0}^{t}f(x,y,\tau)(t-\tau)^{\alpha-1}d\tau,$$
and
$$\Pi y(t)=y_0+\frac{1-\alpha}{AB(\alpha)}g(x,y,t)+\frac{\alpha}{AB(\alpha)\Gamma(\alpha)}\int_{0}^{t}g(x,y,\tau)(t-\tau)^{\alpha-1}d\tau.$$

In page 8, the misprint $\xi=\min\{x_0,x_0\}$ should read $\xi=\min\{u_0,v_0\}$.

In equations (19) and (20),  $\Gamma_1, \Gamma_2$ should simply be $\Pi$, $n\ge 1$.


{\bf Uniqueness}
\begin{Theorem}
	Consider a general Cauchy problem
	\begin{equation*}
	\begin{split}
	^{ABC}_0 D^\alpha_t y(t)&=f(t,y),\;\;\;\;t\in[0, c],[0,c]\subset\mathbb{R}\\
	y(0)&=0,\hspace{1.5cm} y\in\subset([0,c])
	\end{split}
	\end{equation*}
	where $f$ is continuous bounded function on $[0,c]$.
	The above problem has a unique solution if the function $f(t,y(t))$ is Lipschitz with respect to $y$.
\end{Theorem}

{\bf Proof}\\
	To achieve the proof, we define the following mapping
	$$\Pi y=x_0+\frac{1-\alpha}{AB(\alpha)\Gamma(\alpha)}f(t,y(t))+\frac{\alpha}{AB(\alpha)\Gamma(\alpha)}\int_{0}^{t}f(\tau,y(\tau))(t-\tau)^{\alpha-1}d\tau$$
	\begin{equation*}
	\begin{split}
	&\overline{I_c}=[0, c],\\
	&\overline{B_b(y_0)}=[y_0-b, y_0+b]\\
	&C_{c,b}=[0,c]\times \overline{B_b(y_0)}.
	\end{split}
	\end{equation*}
	We adopt the following norm
	$$\|y\|_\infty=\sup_{t\in[0,c]}|y(t)|,\;\;C_{c,b}\longmapsto C_{c,b}$$
	We shall first proof that $(C_{c,b},\|\cdot\|_\infty)$ is Banach space, $C_{c,b}$ is a space of bounded continuous functions. It is well known that the space of bounded functions together with $\|\cdot\|_\infty$ is a Banach space. However, for the sake of readers, we present the proof.
	
	We have that $(I_c,d)$ is a compact metric space, where $d(v,\omega)=\|v-\omega\|$.
	
	For simplicity, let $C_{c,b}=A$, we assume that $\{f_n\}$ is a Cauchy sequence in $A$, thus $\forall\epsilon>0$ there exists $N>0$ such that  $\|f_m-f_n\|_\infty<\epsilon$, $\forall t\in[0,c]$
	$$|f_n(t)-f_m(t)|\le\|f_n-f_m\|_\infty<\epsilon$$
	therefore $\forall t\in[0,c]$ the real numbers $\left\{f_n(t)\right\}$ is a Cauchy sequence. Since $\overline{B_b(y_0)}$ is complete, $\{f_n(t)\}$ is convergent.
	
	Let $$f(t)=\lim\limits_{n\rightarrow\infty}f_n(t).$$
	Here we have $f(t)\in\overline{B_b(y_0)}$. On the other hand $\|f_m\|_\infty<M$, since for every $t\in[0,c]$
	$$|f_n(t)|\le\|f_n\|_\infty\le M.$$
	Now
	$$|f(t)|\le\lim\limits_{n\rightarrow\infty}|f_n(t)|<M,\;\;\forall t\in[0,c]$$
	this shows that $f(t)$ is bounded on $[0,c]$ thus $f\in A$.
	
	$\forall t\in[0,c]$ and $\forall n\ge N$
	$$|f_n(t)-f(t)|=\lim\limits_{m\rightarrow0}|f_n(t)-f_m(t)|<\epsilon$$
	this shows that $\forall n\ge N$, $\|f_n-f\|_\infty<\epsilon$, we have shown that $\{f_n\}$ is convergent in $A$, this proves that $C_{c,b}=A$ with $\|\cdot\|_\infty$ is a Banach space. 	

The above proof can be found in some textbooks. We now have to show that the defined mapping is a contraction under some conditions%
\[
\begin{array}{l}
\left\vert \Pi y-y_{0}\right\vert =\left\vert \frac{1-\alpha }{AB(\alpha )}%
f(t,y(t))+\frac{\alpha }{AB(\alpha )\Gamma \left( \alpha \right) }%
\int\limits_{0}^{t}f(\tau ,y(\tau ))\left( t-\tau \right) ^{\alpha -1}d\tau
\right\vert  \\
\hspace{1.5cm}\leq \frac{1-\alpha }{AB(\alpha )}\left\vert
f(t,y(t))\right\vert +\frac{\alpha }{AB(\alpha )\Gamma \left( \alpha \right)
}\int\limits_{0}^{t}\left\vert f(\tau ,y(\tau ))\right\vert \left( t-\tau
\right) ^{\alpha -1}d\tau
\end{array}%
\]

\[
\begin{array}{l}
\left\Vert \Pi y-y_{0}\right\Vert _{\infty }<\frac{1-\alpha }{AB(\alpha )}%
\sup\limits_{t\in \left[ 0,t\right] }\left\vert f(t,y(t))\right\vert +\frac{%
	\alpha }{AB(\alpha )\Gamma \left( \alpha \right) }\int\limits_{0}^{t}\sup%
\limits_{t\in \left[ 0,t\right] }\left\vert f(\tau ,y(\tau ))\right\vert
\left( t-\tau \right) ^{\alpha -1}d\tau  \\
\hspace{1.5cm}<\frac{1-\alpha }{AB(\alpha )}M+\frac{\alpha M}{AB(\alpha
	)\Gamma \left( \alpha \right) }\int\limits_{0}^{t}\left( t-\tau \right)
^{\alpha -1}d\tau  \\
\hspace{1.5cm}<\frac{1-\alpha }{AB(\alpha )}M+\frac{\alpha M}{AB(\alpha
	)\Gamma \left( \alpha +1\right) }t^{\alpha }<\frac{M}{AB(\alpha )}\left\{
1-\alpha +\frac{c^{\alpha }}{\Gamma \left( \alpha \right) }\right\}
\end{array}%
\]%
$M>\sup\limits_{t\in \left[ 0,t\right] }\left\vert f(t,y(t))\right\vert $
Thus,

\[
\left\Vert \Pi y-x_{0}\right\Vert _{\infty }<\frac{M}{AB(\alpha )}\left\{
1-\alpha +\frac{c^{\alpha }}{\Gamma \left( \alpha \right) }\right\}.
\]%
Here we need%
\[
\left\Vert \Pi y-x_{0}\right\Vert _{\infty }<b
\]%
Thus%
\[
c<\left\{ \left( \frac{AB(\alpha )b}{M}+\alpha -1\right) \Gamma \left(
\alpha \right) \right\} ^{\frac{1}{\alpha }}.
\]%
Next we shall show $\Pi $ is a contraction; $y_{1},y_{2}$ $\in C_{c,b}$

\[
\begin{array}{l}
\left\vert \Pi y_{1}-\Pi y_{2}\right\vert =\left\vert \frac{1-\alpha }{%
	AB(\alpha )}\left( f(t,y_{1})-f(t,y_{2})\right) +\frac{\alpha }{AB(\alpha
	)\Gamma \left( \alpha \right) }\int\limits_{0}^{t}\left( f(\tau
,y_{1})-f(\tau ,y_{2})\right) \left( t-\tau \right) ^{\alpha -1}d\tau
\right\vert  \\
\hspace{1.5cm}\leq \frac{1-\alpha }{AB(\alpha )}\left\vert
f(t,y_{1})-f(t,y_{2})\right\vert +\frac{\alpha }{AB(\alpha )\Gamma \left(
	\alpha \right) }\int\limits_{0}^{t}\left\vert f(\tau ,y_{1})-f(\tau
,y_{2})\right\vert \left( t-\tau \right) ^{\alpha -1}d\tau
\end{array}%
\]

\[
\left\Vert \Pi y_{1}-\Pi y_{2}\right\Vert _{\infty }\leq \frac{1-\alpha }{%
	AB(\alpha )}L\left\Vert y_{1}-y_{2}\right\Vert _{\infty }+\frac{\alpha }{%
	AB(\alpha )\Gamma \left( \alpha \right) }Lt^{\alpha }\left\Vert
y_{1}-y_{2}\right\Vert _{\infty }
\]%
since $f$ is a Lipschitz with respect to $y$

\[
\left\Vert \Pi y_{1}-\Pi y_{2}\right\Vert _{\infty }\leq \frac{L}{AB(\alpha )%
}\left( 1-\alpha +\frac{c^{\alpha }}{\Gamma \left( \alpha \right) }\right)
\left\Vert y_{1}-y_{2}\right\Vert _{\infty }
\]%
To have a contraction, we need
\[
\frac{L}{AB(\alpha )}\left( 1-\alpha +\frac{c^{\alpha }}{\Gamma \left(
	\alpha \right) }\right) <1
\]%
that implies%
\[
c<\left\{ \left( \frac{AB(\alpha )}{L}+\alpha -1\right) \Gamma \left( \alpha
\right) \right\} ^{\frac{1}{\alpha }}.
\]%
Therefore the equation has unique solution if
\[
c<\min \left\{ \left\{ \left( \frac{AB(\alpha )}{L}+\alpha -1\right) \Gamma
\left( \alpha \right) \right\} ^{\frac{1}{\alpha }},\left\{ \left( \frac{%
	AB(\alpha )b}{M}+\alpha -1\right) \Gamma \left( \alpha \right) \right\} ^{%
	\frac{1}{\alpha }}\right\}
\]%
where $AB(\alpha )=1-\alpha +\frac{\alpha }{\Gamma \left( \alpha \right) }$
if $\alpha =1,AB(\alpha )=1.$

{\bf \color{red}

{In Theorem 1},\\
We present a detailed proof of error analysis for approximating the nonlinear function $f(t,x(t),y(t))$ using the Lagrange polynomial interpolation. Finally, we shall consider the case when the fractional order is 1 to show that the function $\Phi(n,\alpha)$ is indeed bounded regardless of the natural number $n\ge1$ chosen.

The error $R^\alpha_n(\zeta)$ is calculated as first integral minus the second integral, that is,
\begin{equation*}
	\begin{split}
		R^\alpha_n(\zeta)&=\frac{\alpha}{AB(\alpha)\Gamma(\alpha)}\int_{0}^{t_{n+1}}(\tau-t_n) (\tau-t_{n-1}) (t_{n+1}-\tau)^{\alpha-1} \left.\frac{\partial^f(\tau,y(\tau))}{\partial\tau^2}\right|_{\tau=\zeta}d\tau \\
		&\;\;\;-\frac{\alpha}{AB(\alpha)\Gamma(\alpha)}\int_{0}^{t_{n}}(\tau-t_n) (\tau-t_{n-1}) (t_{n+1}-\tau)^{\alpha-1} \left.\frac{\partial^f(\tau,y(\tau))}{\partial\tau^2}\right|_{\tau=\zeta}d\tau \\
	\end{split}
\end{equation*}

 In what follows we evaluate the first part of the integral, and also the second part. Bear in mind that $(\tau-t_n)(\tau-t_{n-1})=\tau^2-(t_n+t_{n-1})\tau +t_nt_{n-1}$,  we now evaluate the error as follows:
\begin{equation*}
	\begin{split}
		&\frac{\alpha}{2\Gamma(\alpha)AB(\alpha)}\int_{0}^{n+1}\left[\tau^2-(t_n+t_{n-1})\tau+t_nt_{n-1}\right](t_{n+1}-\tau)^{\alpha-1}\frac{\partial^2}{\partial\tau^2}f(\tau,y(\tau))|_{\tau=\zeta}d\tau\\
		&\le\frac{\alpha}{2\Gamma(\alpha)AB(\alpha)}\int_{0}^{t_{n+1}}\left[\tau^2-(t_n+t_{n-1})\tau +t_nt_{n-1} \right](t_{n+1}-\tau)^{\alpha-1}\sup_{\ell\in[0,\tau]}\left|\frac{\partial^2}{\partial\ell^2}f(\ell,y(\ell))  \right|_{\tau=\zeta}d\tau\\
		&\le \frac{\alpha}{2\Gamma(\alpha)AB(\alpha)} \sup_{t\in[0,t_{n+1}]} \left|\frac{\partial^2}{\partial t^2}f(t,y(t))  \right|_{t=\zeta}
		\int_{0}^{t_{n+1}}\left[\tau^2-(t_n+t_{n-1})\tau +t_nt_{n-1} \right](t_{n+1}-\tau)^{\alpha-1}d\tau\\
		&\le \frac{\sup_{t\in[0,t_{n+1}]} \left|\frac{\partial^2}{\partial t^2}f(t,y(t))  \right|_{t=\zeta}}{2AB(\alpha)(1+\alpha)(2+\alpha)\Gamma(\alpha)}\left\{t_{n+1}^\alpha[(1+\alpha)(2+\alpha)t_nt_{n-1}-(2+\alpha)(t_n+t_{n-1})t_{n+1}+2t^2_{n+1}]  \right\}\\
		&\le \frac{\mathcal{M}}{2AB(\alpha)(1+\alpha)(2+\alpha)\Gamma(\alpha)}t_{n+1}^\alpha\left\{(1+\alpha)(2+\alpha)t_nt_{n-1}-(2+\alpha)(t_n+t_{n-1})t_{n+1}+2t^2_{n+1}]  \right\}
			\end{split}
	\end{equation*}
	We are evaluating the second part of the error error as follows
	\begin{equation*}
		\begin{split}	
		&\frac{\alpha}{2\Gamma(\alpha)AB(\alpha)}\int_{0}^{t_n}[\tau^2-(t_n+t_{n-1})\tau+t_nt_{n-1}]\frac{\partial}{\partial\tau^2}f(\tau,y(\tau))|_{\tau=\zeta}(t_n-\tau)^{\alpha-1}d\tau\\	
		&\frac{\sup_{t\in[0,t_{n}]}\left|\frac{\partial}{\partial t^2}f(t,y(t)) \right|_{t=\zeta}}{2\Gamma(\alpha)AB(\alpha)}
		\int_{0}^{t_n}[\tau^2-(t_n+t_{n-1})\tau+t_nt_{n-1}] (t_n-\tau)^{\alpha-1}d\tau\\
		&\le \frac{\mathcal{M}\alpha t_n^{\alpha+1}}{2AB(\alpha)(1+\alpha)(2+\alpha)\Gamma(\alpha)}\{(2+\alpha)t_{n-1}-t_n \}
	\end{split}
\end{equation*}

Therefore,
\begin{equation*}
	\begin{split}
		&|R^\alpha_n(\zeta)|\le \frac{\mathcal{M}}{2(1+\alpha)(2+\alpha)\Gamma(\alpha)AB(\alpha)}\left|  \left\{2\alpha t^{2+\alpha}_n+t_{n+1}^{1+\alpha}[2t_{n+1}-(2+\alpha)t_n] +(2+\alpha)t_{n-1}\right.\right.\\
		&\hspace{2cm}\left.\left.\left[-2\alpha  t^{\alpha+1}_n+((1+\alpha)t_n-t_{n+1} )t^\alpha_{n+1} \right]   \right\}\right|
	\end{split}
\end{equation*}

\begin{equation*}
	\begin{split}
		&|R^\alpha_n(\zeta)|\le \frac{\mathcal{M}\Delta t^{\alpha+2}}{2(1+\alpha)(2+\alpha)\Gamma(\alpha)AB(\alpha)}
		\left|\left\{2\alpha n^{2+\alpha}+(n+1)^{\alpha+1}(2(n+1)\right.\right.\\
		&\hspace{1.5cm}\left.\left.-(2+\alpha)n)+(2+\alpha)(n-1)[-2\alpha n^{\alpha+1}((1+\alpha)n-(n+1))(n+1)^\alpha] \right\}\right|
	\end{split}
\end{equation*}

\begin{equation*}
	\begin{split}
		&|R^\alpha_n(\zeta)|\le \frac{\mathcal{M}\Delta t^{\alpha+2}_\alpha}{2\Gamma(\alpha+3)}
		\left|\left\{2\alpha n^{2+\alpha}+(n+1)^{\alpha+2}(2-\alpha n)\right.\right.\\
		&\hspace{1.5cm}\left.\left.+(2+\alpha)(n-1)[-2\alpha n^{\alpha+1}+(\alpha n-1)(n+1)^\alpha] \right\}\right|
	\end{split}
\end{equation*}

$$|R^\alpha_n(\zeta)|<\frac{\mathcal{M}\Delta t^{\alpha}}{2\Gamma(\alpha+3)}\Phi(n,\alpha)$$
where
\begin{equation*}
	\begin{split}
		&|\Phi(n,\alpha)|\\
		&\;\;\;=\left|\alpha\left\{2\alpha n^{2+\alpha}+(n+1)^{\alpha+2}(2-\alpha n) +(2+\alpha)(n-1)[-2\alpha n^{\alpha+1}+(\alpha n-1)(n+1)^\alpha] \right\}\right|.
	\end{split}
\end{equation*}


 Indeed the error depends on $\alpha$, therefore for a given $\alpha$ we get an associate error, which is  different for the case of classical  integral as the order is fixed to $\alpha=1$, the error depends only on $\Delta h$. It should be noted that a fractional integral has memory which is depicted by $ (t_{n+1}-\tau)^{\alpha-1}$.
 In what follows, we shall give some  examples}

\section*{Example 1}
We consider the following Cauchy problem with ABC derivative
\begin{equation}\label{ex1}
\begin{split}
^{ABC}_0 D^\alpha_t y(t)&=t^\beta,\;\;\;\;[0, 10],\\
y(0)&=0
\end{split}
\end{equation}
and
$$y(t)=\frac{1-\alpha}{AB(\alpha)}t^\beta+\frac{\alpha t^{\alpha+\beta}\Gamma(\beta+1)}{AB(\alpha)\Gamma(\alpha+\beta+1)}.$$

Although the above problem has the exact solution, we will verify the convolution
$$\Pi y(t)=\frac{(1-\alpha)}{AB(\alpha)}f(t,y(t))+\frac{\alpha}{AB(\alpha)\Gamma(\alpha)}\int_{0}^{t}f(\tau,y(\tau))(t-\tau)^{\alpha-1}d\tau$$
\begin{equation*}
\begin{split}
f(t,y(t))&=t^\beta\\
|f(t,y(t))|&=|t^\beta|\le 10^\beta
\end{split}
\end{equation*}

\begin{equation*}
\begin{split}
|\Pi y(t)-y(0)|&<\frac{(1-\alpha)10^\beta}{AB(\alpha)}+\frac{\alpha}{AB(\alpha)\Gamma(\alpha)}\int_{0}^{t}t^\beta(t-\tau)^{\alpha-1}d\tau\\
&<\frac{(1-\alpha)10^\beta}{AB(\alpha)}+\frac{\alpha 10^\beta}{AB(\alpha)\Gamma(\alpha)} \frac{t^\alpha\Gamma(\alpha)}{\Gamma(\alpha+1)}\\
&<\frac{(1-\alpha) 10^\beta}{AB(\alpha)}+\frac{\alpha 10^{\beta+\alpha}}{AB(\alpha)\Gamma(\alpha+1)}<10
\end{split}
\end{equation*}
Since $0<\alpha\le 1$, then
$$|\Pi y_1-\Pi y_2|=0<\frac{1}{1000}|y_1-y_2|$$
$k=\frac{1}{1000}<1$
which is a contraction. \\
We can conclude that the solution exists and unique.

\begin{figure}[!h]
	\centering
	\begin{tabular}{c}
		\begin{minipage}{300pt}
			\includegraphics[width=300pt]{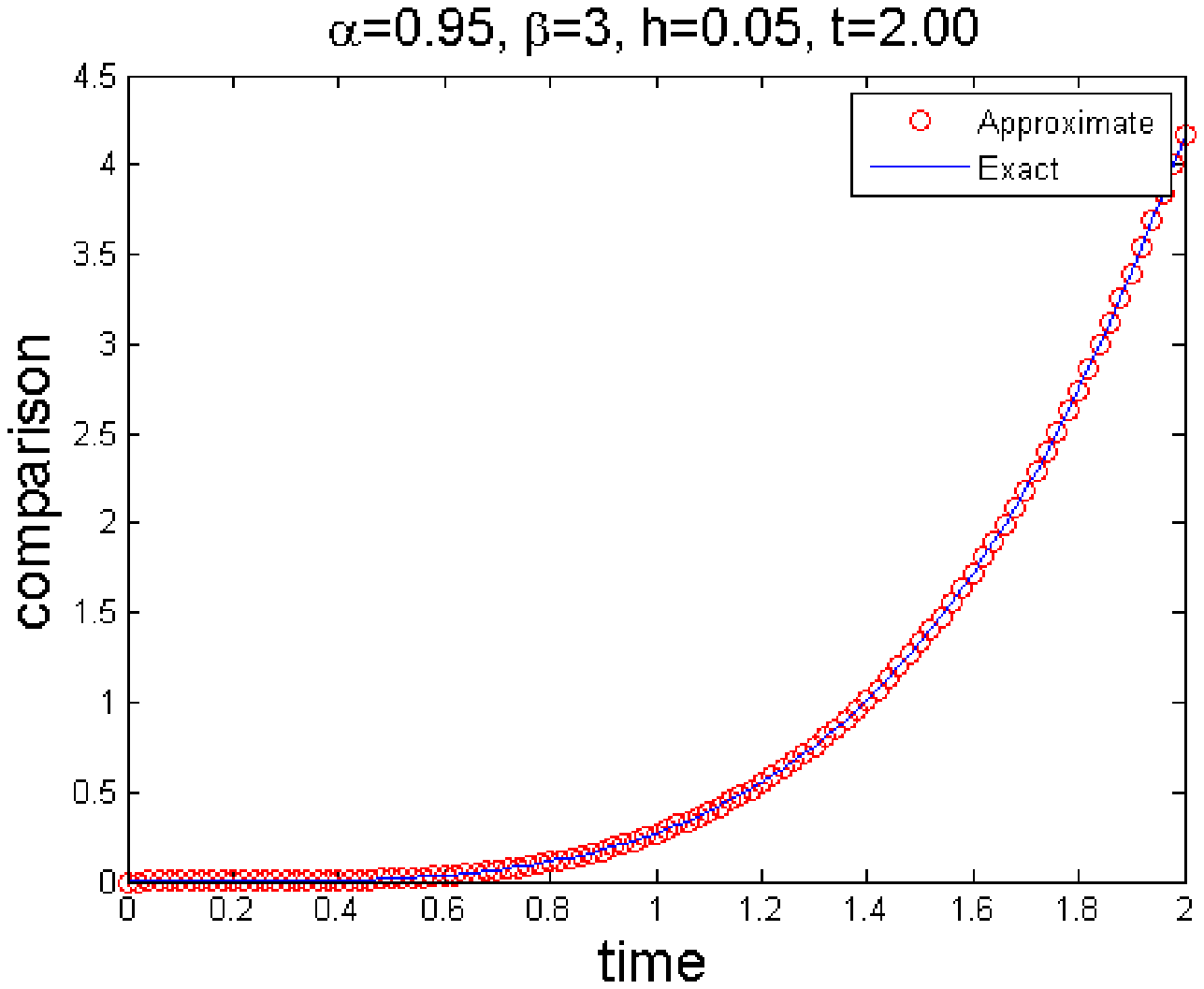}
		\end{minipage}\\\\
		\begin{minipage}{300pt}
			\includegraphics[width=300pt]{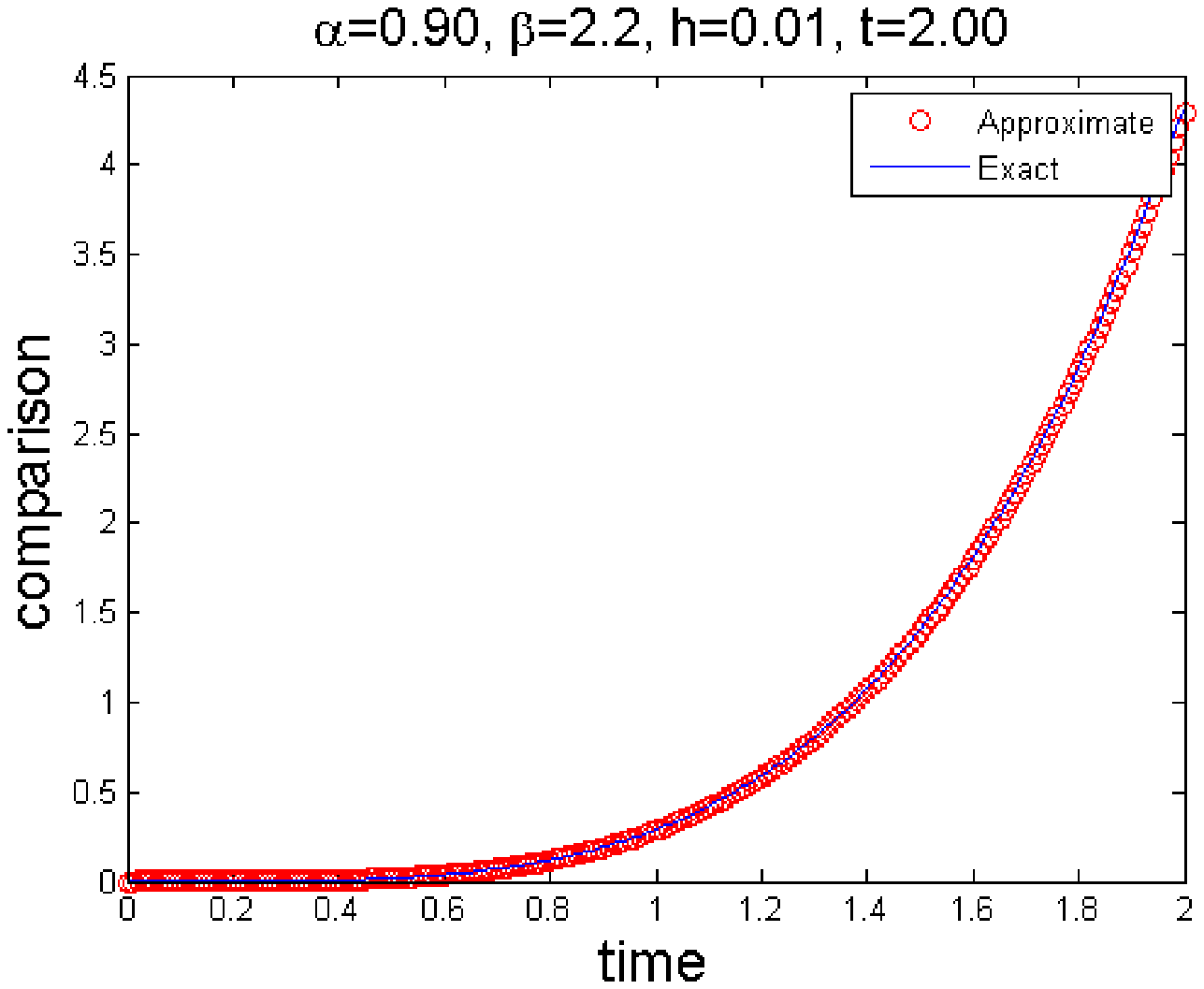}
		\end{minipage}
	\end{tabular}
	\caption{Numerical solutions for test problem 1 for different values of $\alpha,\beta$ and $h$ with $t=2$.}\label{F1}
\end{figure}

The numerical results showing comparison between the exact and approximate solutions for different instances of fractional index $\alpha$ and $\beta$ for $t=2$ is displayed in Figure \ref{F1}.

\section*{Example 2}
We consider the following
\begin{equation}\label{ex2}
\begin{split}
^{ABC}_0 D^\alpha_t y(t)&=\frac{t}{1000}y(t),\;\;\;\;[0, 10],\\
y(0)&=0
\end{split}
\end{equation}
$$\Pi y(t)=\frac{1-\alpha}{AB(\alpha)}\frac{ty(t)}{1000}+\frac{\alpha}{AB(\alpha)\Gamma(\alpha)}\int_{0}^{t}\tau y(\tau)(t-\tau)^{\alpha-1}d\tau$$
$$|\Pi y_1-\Pi y_2|=\left|\frac{t}{1000}(y_1-y_2)(1-\alpha)+\frac{\alpha}{AB(\alpha)\Gamma(\alpha)}\int_{0}^{t}\frac{\tau(t-\tau)^{\alpha-1}}{1000}(y_1-y_2)d\tau \right|$$
Nothing that $f(t,y(t))=\frac{ty(t)}{1000}$
\begin{equation*}
\begin{split}
|f(t,y_1)-f(t,y_2)|&=\frac{t}{1000}|y_1-y_2|\\
&\le\frac{t}{100}|y_1-y_2|
\end{split}
\end{equation*}
where $f$ is Lipschitz with $k=1/100$
\begin{equation*}
\begin{split}
\|\Pi y_1-\Pi y_2\|_\infty&\le\frac{t}{1000}\|y_1-y_2\|_\infty\frac{(1-\alpha)}{AB(\alpha)}+\frac{\alpha}{AB(\alpha)\Gamma(\alpha)}\|y_1-y_2\|_\infty \frac{1}{100}\int_{0}^{t}(t-\tau)^{\alpha-1}d\tau\\
|\Pi y_1-\Pi y_2\|_\infty&\le\frac{(1-\alpha)}{100 AB(\alpha)}\|y_1-y_2\|_\infty+\frac{\alpha}{100AB(\alpha)} \|y_1-y_2\|_\infty\frac{10^\alpha}{\Gamma(\alpha+1)}\\
&<\frac{1}{100}\left\{\frac{1-\alpha}{AB(\alpha)}+\frac{\alpha 10^\alpha}{AB(\alpha)\Gamma(\alpha+1)}\right\} \|y_1-y_2\|_\infty\\
&\le k\|y_1-y_2\|_\infty
\end{split}
\end{equation*}
$k<1$, which is a contraction.
\begin{equation}
\begin{split}
\Pi y-y(0)&=\frac{1-\alpha}{AB(\alpha)}\frac{ty}{1000}+\frac{\alpha}{AB(\alpha)\Gamma(\alpha)}\int_{0}^{t}\tau y(\tau)\frac{(t-\tau)^{\alpha-1}}{1000}d\tau\\
\|\Pi y-y(0)\|_\infty&\le\frac{(1-\alpha)}{AB(\alpha)100}\|y\|_\infty+\frac{\alpha\|y\|_\infty 10^\alpha}{AB(\alpha)1000\Gamma(\alpha+1)}\\
&\le\left(\frac{1-\alpha}{AB(\alpha)}+\frac{10^{\alpha-1}}{\Gamma(\alpha+1)AB(\alpha)}\right)\frac{\|y\|_\infty}{100}
\end{split}
\end{equation}
Hence, we can conclude that the solution exists and unique.

\subsection*{ Error for classical- vs fractional-derivatives in Caputo sense}
In this section, we present a derivation of the error analysis for the case of $\alpha =1$.  This will help the reader to see the difference between the Cauchy problem with fractional derivative and that with classical derivative. It is worth noting the error for the Caputo case is the same as the Atangana-Baleanu derivative case since the fractional integral associated with the Atangana-Baleanu operator is the same with an ${\alpha}/{AB(\alpha)}$ coefficient.
	

Assuming that a function $f(t,y(t))$ is two time differentiable, such that $\frac{\partial^2f(t,y(t))}{\partial t^2}$ is bounded.

\begin{equation}\label{eq1}
\begin{split}
\frac{dy(t)}{dt}=f(t,y(t)),\;\;[0,T]\\
y(0)=y_0
\end{split}
\end{equation}
From (\ref{eq1}) we apply the integral to get
$$y(t)-y(0)=\int_{0}^{t}f(\tau,y(\tau))d\tau,$$
$t_{n+1}=(n+1)\Delta t$ and $t_n=n\Delta t$
So that
\begin{equation}\label{eq2}
y(t_{n+1})-y(0)=\int_{0}^{t}f(\tau,y(\tau))d\tau
\end{equation}
and
\begin{equation}\label{eq3}
y(t_n)-y(0)=\int_{0}^{t_n}f(\tau,y(\tau))d\tau
\end{equation}
Taking (\ref{eq2})-(\ref{eq3}) yields
\begin{equation}
y(t_{n+1})-y(t_n)=\int_{t_n}^{t_{n+1}}f(\tau,y(\tau))d\tau
\end{equation}
where $f(\tau,y(\tau))$is approximate within $[t_n,t_{n+1}]$ with Lagrange polynomial
\begin{equation}
\begin{split}
f(\tau,y(\tau))=P_n(\tau)&=\frac{\tau-t_{n-1}}{\Delta t}f(t_n,y(t_n))-\frac{\tau-t_{n}}{\Delta t}f(t_{n-1},y(t_{n-1}))\\
&\;\;\;\;+\left.\frac{\partial^2f(\tau,y(\tau))}{\partial \tau^2}\right|_{\tau=\zeta} \frac{(\tau-t_{n-1})}{2!}.
\end{split}
\end{equation}
The error is given as
$$R_n(\zeta)=\int_{t_n}^{t_{n+1}}\left.\frac{\partial^2f(\tau,y(\tau))}{\partial \tau^2}\right|_{\tau=\zeta}\frac{(\tau-t_n)(\tau-t_{n-1})}{2!}d\tau$$

\begin{equation*}
\begin{split}
|R_n(\zeta)|&=\left|\int_{t_n}^{t_{n+1}} \left.\frac{\partial^2f(\tau,y(\tau))}{\partial \tau^2}  \right|_{\tau=\zeta} \frac{(\tau-t_n)(\tau-t_{n-1})}{2!}d\tau\right|\\
&\le \int_{t_n}^{t_{n+1}} \left|\left.\frac{\partial^2f(\tau,y(\tau))}{\partial \tau^2}  \right|_{\tau=\zeta}\right| \frac{(\tau-t_n)(\tau-t_{n-1})}{2}d\tau\\
&<\int_{t_n}^{t_{n-1}}\sup_{\ell\in[t_n,\tau]}
\left|\frac{\partial^2f(\tau,y(\tau))}{\partial \tau^2}  \right|_{\tau=\zeta} \frac{(\tau-t_n)(\tau-t_{n-1})}{2}d\tau
\end{split}
\end{equation*}
Since we assume that $\frac{\partial^2}{\partial t^2}$ is bounded
\begin{equation*}
\begin{split}
|R_n(\zeta)|&<\sup_{t\in[t_n,t_{n+1}]}  \left|\frac{\partial^2f(t,y(t))}{\partial \tau^2}  \right|_{t=\zeta} \int_{t_n}^{t_{n-1}} \frac{(\tau-t_n)(\tau-t_{n-1})}{2!}d\tau\\
&<\frac{M}{2}\left[\frac{\tau^3}{3}-(t_{n-1}+t_n)\frac{\tau^2}{2}+t_{n-1}t_n\tau \right]_{t_n}^{t{n-1}}\\
&<\frac{M}{2}\left\{\frac{t_{n+1}^3}{3}-\frac{t_n^3}{3}-(t_{n-1}+t_n)\left[\frac{t_{n+1}^2}{2}-\frac{t^2_n}{2} \right]+t_{n-1}t_n[t_{n+1}-t_n]    \right\}\\
&<\frac{M}{2}\Delta t^3 \left\{\frac{(n+1)^3-n^3}{3}-(2n-1)\left(\frac{(n+1)^2}{2}-\frac{n^2}{2} \right)+(n-1)n \right\}\\
&<\frac{M}{2}\Delta t^3\left\{\frac{3n^2+3n+1}{3}-\frac{(4n^2-1)}{2}+(n-1)n \right\}\\
&<\frac{M}{2}\Delta t^3\left\{\frac{1}{3}+\frac{1}{2}\right\}\\
&<\frac{5}{12}M\Delta t^3.
\end{split}
\end{equation*}
The above error depends on $\Delta t$ only because of integer order. So the numerical scheme for this is
$$y_{n+1}=y_n+\frac{3}{2}\Delta tf(t_n,y_n)-\frac{\Delta t}{2}f(t_{n-1},y_{n-1}).$$
This method is well-known to be accurate and is called Adams-Bashforth scheme.\\
What is the difference with the fractional case?
\begin{equation*}
\begin{split}
y(t_{n+1})-y(t_n)&=\frac{1}{\Gamma(\alpha)}\int_{0}^{t_{n+1}}(t_{n+1}-\tau)^{\alpha-1} f(\tau,y(\tau))d\tau \\
&\;\;-\frac{1}{\Gamma(\alpha)}\int_{0}^{t_n}(t_n-\tau)^{\alpha-1}f(\tau,y(\tau))d\tau.
\end{split}
\end{equation*}
We have new component which is the power law contribution $(t_{n+1}-\tau)^{\alpha-1}$ and $(t_{n}-\tau)^{\alpha-1}$, by approximating $f(t,y(t))$

 \begin{equation*}
 \begin{split}
 y(t_{n+1})&=y(t_n)+\left[\int_{0}^{t_{n+1}}(t_{n+1}-\tau)^{\alpha-1}\left\{\frac{(\tau-t_{n-1})}{\Delta t} f(t_n,y(t_n))\right.\right.\\
 &\;\;\left.\left.-\frac{\tau-t_n}{\Delta t} f(t_{n-1},y(t_{n-1}))  \right\}d\tau-\int_{0}^{t_n}(\tau-t_n)(\tau-t_{n-1})(t_n-\tau)^{\alpha-1}d\tau \right]\frac{1}{\Gamma(\alpha)} \\
 &\;\;+\left\{ \int_{0}^{t_{n+1}}(t_{n+1}-\tau)^{\alpha-1}\frac{\partial^2f(\tau,y(\tau))}{\partial \tau^2}|_{\tau=\zeta}(\tau-t_n)(\tau-t_{n-1}) d\tau\right.\\
 &\;\;\left. - \int_{0}^{t_{n}}(t_{n}-\tau)^{\alpha-1}\frac{\partial^2f(\tau,y(\tau))}{\partial \tau^2}|_{\tau=\zeta}(\tau-t_n) \frac{\tau-t_{n-1}}{2!}d\tau \right\}
 \end{split}
 \end{equation*}
 These $(t_{n+1}-\tau)^{\alpha-1}$ and $(t_{n}-\tau)^{\alpha-1}$ are nonlinear if $\alpha=1$, indeed
 $$ y(t_{n+1})-y(t_n)=\int_{t_n}^{t_{n+1}}(\tau-t_n)(\tau-t_{n-1}) \frac{\partial^2f(\tau,y(\tau))}{\partial \tau^2} |_{t=\tau}d\tau+\int_{t_n}^{t_{n+1}}f(\tau,y(\tau))d\tau $$
 Then, we have
 $$|R_n(\zeta)|<\frac{5M\Delta t}{12}.$$
 But if $\alpha\ne 1$, $0<\alpha<1$ then
 \begin{equation*}
 \begin{split}
 R_n(\zeta)&=\frac{1}{\Gamma(\alpha)}\left\{\int_{0}^{t_{n+1}}(t_{n+1}-\tau)^{\alpha-1}\frac{(\tau-t_n)(\tau-t_{n-1})}{2!}\frac{\partial^2f}{\partial \tau^2}(\tau,y(\tau))|_{\tau=\zeta}d\tau  \right.\\
 &\;\;\left.- \int_{0}^{t_{n}}(t_{n+1}-\tau)^{\alpha-1}\frac{(\tau-t_n)(\tau-t_{n-1})}{2!}\frac{\partial^2f}{\partial \tau^2}(\tau,y(\tau))|_{\tau=\zeta}d\tau  \right\}\\
 |R_n(\zeta)|&<\frac{Mh^{\alpha+2}}{2!}\Phi(n,\alpha)
 \end{split}
 \end{equation*}
 Indeed, if $\alpha=1$, then $\Phi(n,\alpha)=\frac{5}{12}$.

 {\bf \color{red} In Figure \ref{F2} below, we have presented the plots of $\Phi(n,\alpha)\cdot\frac{h^{\alpha+2}}{2}$, first as a function of $n$ for different values of alpha from $0.1$ to $1$, then as function of $n$ and $\alpha$. The results show that when $\Phi(n, \alpha)<\alpha(n,1)$. }

\begin{figure}[!h]
	\centering
	\begin{tabular}{c}
		A
		\begin{minipage}{300pt}
			\includegraphics[width=300pt]{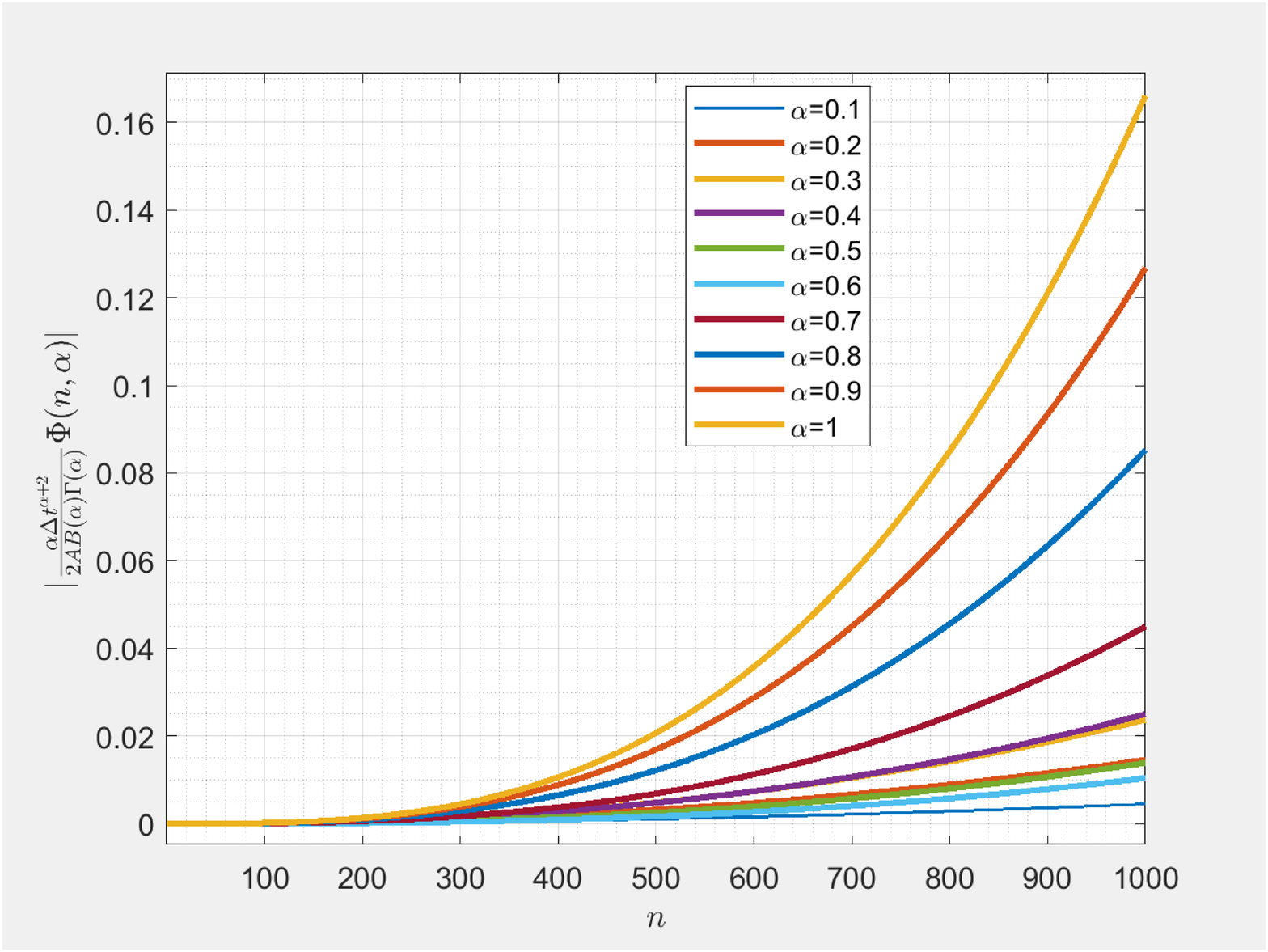}
		\end{minipage}\\\\
	B
		\begin{minipage}{300pt}
			\includegraphics[width=300pt]{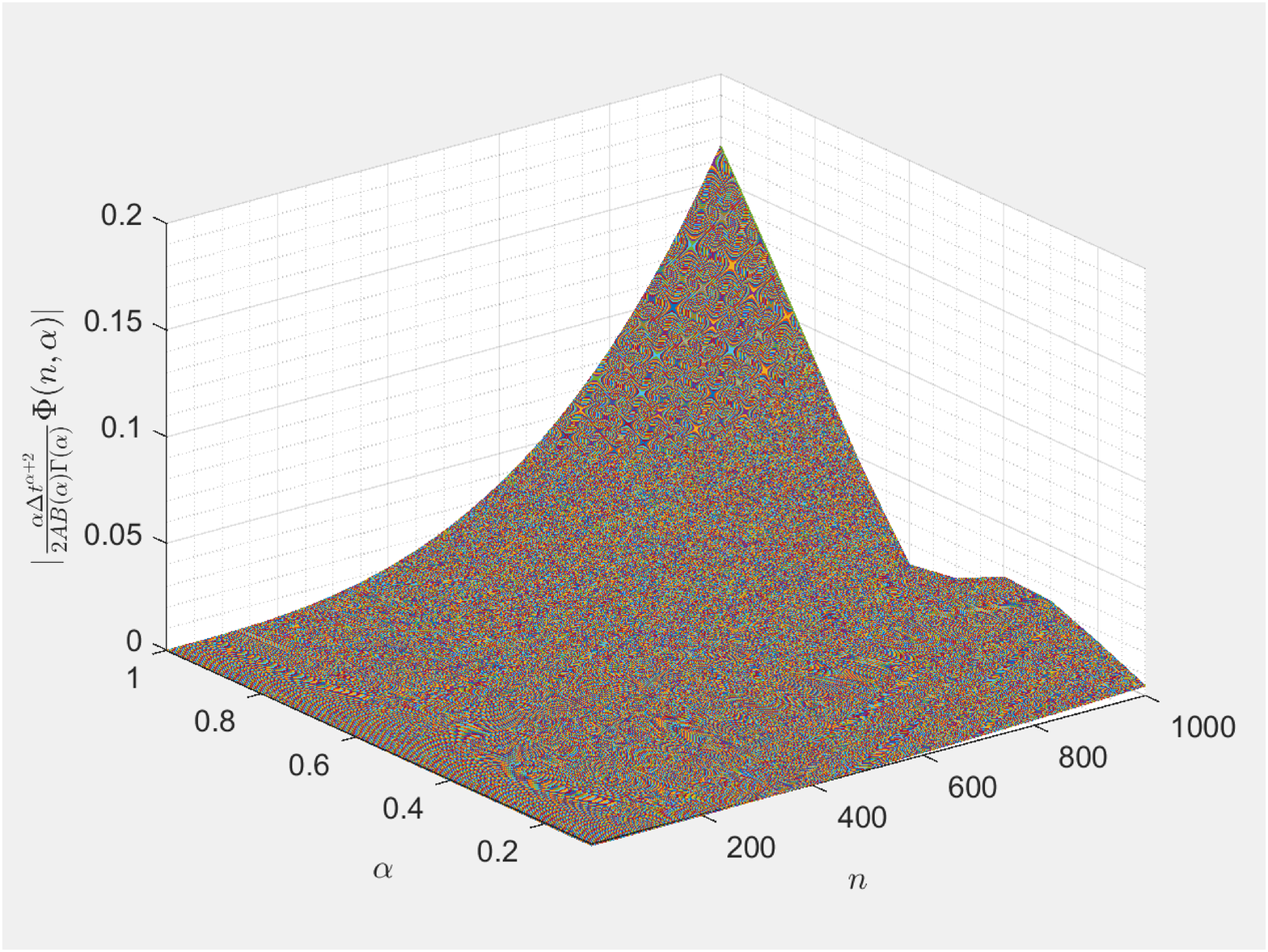}
		\end{minipage}
	\end{tabular}
	\caption{{\bf \color{red}A: shows plot of $\Phi(n,\alpha)$ as a function of  $n$ for different values of $\alpha$. B: denotes $\Phi(n,\alpha)$ as a function of $n$ and $\alpha$.}}\label{F2}
\end{figure}

{\bf \color{blue}
	Alternatively, one can proceed as follows:
\begin{equation}
	\begin{split}
		u(t_{n+1})-u(t_n)&=\frac{1-\alpha}{AB(\alpha)}\left[ f(t_n,u(t_n))- f(t_{n-1},u(t_{n-1}))  \right]\\
		&+ \frac{\alpha}{AB(\alpha)\Gamma(\alpha)}\sum_{j=1}^{n} \int_{t_j}^{t_{j+1}} f(\tau,u(\tau)) (t_{n+1}-\tau)^{\alpha-1}d\tau\\
		&-\frac{\alpha}{AB(\alpha)\Gamma(\alpha)} \sum_{j=1}^{n-1}\int_{t_j}^{t_{j+1}} f(\tau,u(\tau)) (t_n-\tau)^{\alpha-1}d\tau
	\end{split}
\end{equation}
such that in $[t_j,t_{j+1}]$ the function $f(\tau,u(\tau))$ can be interpolated using the Lagrange polynomial interpolation for the first integral. Also, the function can be interpolated within $[t_k,t_{k+1}]$ for the second integral. }

\section*{Acknowledgment}
The authors are grateful to all of the anonymous reviewers and the Editor-in-Chief for their professional support and valuable comments.

\end{document}